\documentclass[11pt,a4paper]{article}
\usepackage{srcltx}
\usepackage{amsmath}
\usepackage{multirow}
\usepackage{slashbox}
\def\Size{\hbox{\rm Size}}
\def\mes{\hbox{\rm mes}}
\usepackage[bf,normalsize,tableposition=top]{caption}
\usepackage{subfig}
\usepackage{amsfonts}
\usepackage{color}
\usepackage{epsfig}
\usepackage{amssymb}
\usepackage{graphicx}
\usepackage{algorithm}
\usepackage[colorlinks=true]{hyperref}

\newcommand{\cP}{{\cal P}}

\def\H{{\cal H}}

\def\T{{\cal T}}

\def\cR{{\cal R}}

\def\Opt{{\mathop{\hbox{\rm Opt}}}}

\newcommand{\bbr}{{\mathbb{R}}}
\newcommand{\bbz}{{\mathbb{Z}}}

\newcommand{\cN}{I\!\! N}

\def\M{{\cal M}}
\def\A{{\cal A}}
\def\cP{{\cal P}}

\newcommand{\Z}{{\cal Z}}

\newcommand{\F}{{\cal F}}
\newcommand{\E}{{\cal E}}
\newcommand{\C}{{\cal C}}

  \newcommand{\X}{{\cal X}}

\newcommand{\half}{ \mbox{\small$\frac{1}{2}$}}

\newcommand{\four}{ \mbox{\small$\frac{1}{4}$}}

\def\inter{{\hbox{\rm int}\,}}

\def\Prob{{\hbox{\rm Prob}}}

\newcommand{\be}{\begin{eqnarray}}
\newcommand{\ee}[1]{\label{#1}\end{eqnarray}}

\newcommand{\ese}{\end{eqnarray*}}
\newcommand{\bse}{\begin{eqnarray*}}
\newcommand{\rf}[1]{~(\ref{#1})}

\newtheorem{proposition}{Proposition}
\newtheorem{corollary}{Corollary}
\newtheorem{theorem}{Theorem}

\newtheorem{remark}{Remark}

\def\Risk{\hbox{\rm Risk}}

\def\cN{{\cal N}}
\def\J{{\cal J}}
\def\qed{\hfill$\square$}

\def\O{{\cal O}}

\title{On sequential hypotheses testing via convex optimization}
\author{
Anatoli Juditsky
\thanks{LJK, Universit\'e Grenoble Alpes, B.P. 53, 38041 Grenoble Cedex 9, France,	
{\tt anatoli.juditsky@imag.fr}}
\and Arkadi Nemirovski
\thanks{Georgia Institute
 of Technology, Atlanta, Georgia
30332, USA, {\tt nemirovs@isye.gatech.edu}\newline
Research of the first author was supported by the CNRS-Mastodons project GARGANTUA,
and the LabEx PERSYVAL-Lab (ANR-11-LABX-0025). Research of
the second author was supported by NSF grants CMMI-1232623, CMMI-1262063, CCF-1415498.}
}
\begin{document}
\maketitle
\begin{abstract}
We propose a new approach to sequential testing which is an adaptive (on-line) extension of the (off-line) framework developed in \cite{GJuN2014}. It relies upon testing {of} pairs of hypotheses in the case where each hypothesis states that the vector of parameters underlying the distribution of observations belongs to a convex set. The nearly optimal under appropriate conditions test is yielded by a solution to an efficiently solvable convex optimization
problem. The proposed methodology can be seen as a computationally friendly reformulation of the classical sequential testing.
\end{abstract}
\section{Introduction}
{
Let $\omega_i\in\Omega$, $i=1,2,...$  be a sequence of independent and identically distributed (iid) observations with a common density w.r.t. a given measure $P$ on the observation space $\Omega$, known to belong to a given parametric family $\{p_\mu(\cdot):\mu\in\M\}$. Assume that we are given a (finite) number $I$ of subsets $X^i$, $1\leq i\leq I$, of the parameter space $\M$; these sets define hypotheses $H_i$ on the density $p(\cdot)$ underlying our observation, $H_j$ stating that
the parameter of this density belongs to $H_j$: $p(\cdot)=p_\mu(\cdot)$ with $\mu\in X^i$. We are interested in the problem of multiple testing of  composite hypotheses: our goal is, given observations $\omega_i$ to decide on the hypotheses $H_1,...,H_I$. When the size $K$ of the observation $\omega^K=(\omega_1,...,\omega_K)$ is fixed, the performance of a test $T^K(\cdot)$ -- a measurable mapping of $\Omega^K$   into $\{1,...,I\}$ -- can be quantified by the maximal  error probability of rejecting the true hypothesis:
\be
\Risk(T^K)=\max_{1\leq i\leq I} \sup_{\mu\in X^i}  \Prob_{\mu,K}\{\omega^K:\;T^K(\omega^K)\neq i\}
\ee{risk1}
where $\Prob_{\mu,K}$ is the probability w.r.t. the distribution of the observation $\omega^K$ corresponding to the ``true'' parameter $\mu\in \M$. Then the problem of optimal testing can be approached through minimization of the risk over the class of tests.
Yet, in many practical applications, especially in those where observations come at a price, such approach may be too conservative.
The sequential approach to testing, introduced in the pioneering papers of Barnard \cite{barnard46} and Wald \cite{wald45,ww48} is now a reach area of statistical theory and offers many strong results (see also \cite{Chernoff1972,Ghosh1991,Bakeman1997} for references or \cite{lai01} for a recent review). For the problem of interest this approach can be summarized as follows: given an upper bound $\epsilon$ of the risk the problem is reexamined sequentially, when new observations are available. The test terminates when either a decision with the risk $\leq \epsilon$ is possible, or when maximal allowed observation count $K$ is reached (in which case the test produces no decision). Usually, the sequential test is based on the {\em Generalized Likelihood} or {\em Mixture Likelihood} statistics, and its performance is evaluated by its asymptotical, as $\epsilon\to+0$ and $K\to\infty$. The approach we pursuit in this paper is of completely different spirit. It originates from \cite{JN2009} and was applied to the testing problem in \cite{GJuN2014} in the case where the size of the observation sample is fixed. The main ``building block'' of this approach is
a construction, based on Convex Programming (and thus computationally efficient) allowing,
under appropriate assumptions, to build a provably nearly optimal test for deciding, given observation $\omega^K$, between a
pair of composite hypotheses of the sort $H_1:\,\mu\in X$ and $H_2:\,\mu\in Y$ where $X$ and $Y$ are convex compact subsets of $\M$.
This approach is applicable in several important situations, namely where (a) $p_\mu$ is Gaussian density on $\bbr^n$ with expectation $\mu$ and fixed covariance matrix,    (b) $p_\mu$ is the distribution of the Poisson vector in $\bbr^m$ with independent components with parameters $\mu_i$, and (c) is a distribution of a discrete random variable taking values in $\{1,...,m\}$ with $\mu\in \bbr_+$ being the vector of probabilities of the corresponding distribution. As a compensation for rather restrictive assumptions on the families of densities $p_\mu$, the approach in question is extremely permissive as far the structure of the sets $X$ and $Y$ is concerned: what we require (apart from compactness and convexity) is the efficient tractability of $X$ and $Y$.
While the analytical expressions for the test characteristics  are not available in the proposed approach,
rather detailed information about performance guarantees of the resulting decisions can be obtained by efficient situation-oriented computation.
In \cite{GJuN2014} we introduced a calculus of pairwise tests to design nearly optimal testing procedures in the situation where the sets $X^i$ corresponding to different hypotheses $H^i$ are unions of (not too large number of) convex sets with the risk of the test defined as in \rf{risk1}  (in the minimax setting).\footnote{It should be mentioned that what we call below ``simple
tests'' were used to test composite hypotheses represented by convex sets in the white noise model (a) in  \cite{Burnashev1979,Burnashev1982,Ingster2002} and in the distribution model (c) in  \cite{Lecam1973,Lecam1975,Birge1980,Birge1982,Birge1983M,Lecam1986}.} What follows can be seen as an adaptive version of the testing procedure from \cite{GJuN2014} -- when the ``true distribution'' of observations corresponds to the value $\mu$ of the parameter which is  ``deeply inside'' of some $X^i$, the correct decision ($H_i$ is accepted) can be taken much faster (using a smaller observation sample) than in the case of $\mu$ close to some ``wrong'' $H^i$'s. On the other hand, it is nothing but a ``computationally friendly'' application of the classical Sequential Analysis methodology to the problem in question.
\par
The paper is organized as follows. In section \ref{sec:oldies} we give a summary of the results of \cite{GJuN2014} on {\sl hypothesis testing in good observation schemes} and discuss the construction of quasi-optimal ``off-line'' test. Then in section \ref{sec:seq1} we introduce the sequential problem setting and construct a generic test aggregation procedure which allows to reduce multiple testing to pairwise testing. Finally, in section \ref{sec:impl1} we consider in detail a particular implementation of the proposed approach and present some very preliminary simulation results.
The proofs missing in the main body of the paper are postponed to the appendices.}
\section{Preliminaries}\label{sec:oldies}
 What follows is the summary of the approach of \cite{JN2009} as applied to hypotheses testing; for detailed presentation of the constructions and results of this section, same as for the related proofs, see \cite{GJuN2014}.
\subsection{Good observation schemes}\label{sec:good_os}
We start with introducing {\sl good observation schemes}, those with which we intend to work. Recall that we are interested to make inferences from a random observation $\omega$ taking values in a given {\sl observation space} $\Omega$ and obeying probability density w.r.t. a given measure $P$ on $\Omega$; this density is known to belong to a given parametric family $\{p_\mu(\cdot):\mu\in\M\}$ of probability densities, taken w.r.t. $P$, on $\Omega$.
We intend to work under the following assumptions on our ``observation environment:''
\begin{quote}
\begin{enumerate}
\item $\M\subset\bbr^m$ is a convex set which coincides with its relative interior;
\item $\Omega$ is a Polish (i.e., separable complete metric) space equipped with a Borel $\sigma$-additive $\sigma$-finite measure $P$, $\hbox{\rm supp}(P)=\Omega$, and { distributions $P_\mu\in\cP$ possess densities $p_\mu(\omega)$ w.r.t. $P$. We assume that}
    \begin{itemize}
    \item $p_\mu(\omega)$ is continuous in $\mu\in\M$, $\omega\in\Omega$ and is positive;
    \item the densities $p_\mu(\cdot)$ are ``locally uniformly summable:'' for every compact set $M\subset\M$,  there exists a Borel function $p^M(\cdot)$ on $\Omega$ such that
    $\int_\Omega p^M(\omega)P(d\omega)<\infty$ and $p_\mu(\omega)\leq p^M(\omega)$ for all $\mu\in M$, $\omega\in\Omega$;
    \end{itemize}
\item We are given a finite-dimensional linear space $\F$ of continuous functions on $\Omega$ containing constants such that $\ln(p_\mu(\cdot)/p_\nu(\cdot))\in\F$ whenever $\mu,\nu\in\M$.
    \par Note that the latter assumption implies that  distributions $P_\mu,\;\mu\in\M$, belong to an exponential family.
\item
For every $\phi\in \F$, the function
$
F_\phi(\mu)=\ln\left(\int_\Omega\exp\{\phi(\omega)\}p_\mu(\omega)P(d\omega)\right)
$
is well defined and concave in $\mu\in \M$.
 \end{enumerate}
 \end{quote}
In the just described situation, where assumptions 1-4 hold, we refer to the collection $\O=((\Omega,P),\{p_\mu(\cdot):\mu\in\M\},\F)$  as  {\em good observation scheme} (o.s.).
\subsubsection{Basic examples}\label{sect:examples}
Basic examples of good o.s.'s are as follows:
\paragraph{Gaussian o.s.} Here $\Omega=\bbr^m$, $P$ is the Lebesque measure on $\Omega$, $\M=\bbr^m$ and $p_\mu(\omega)=\cN(\mu,I_m)$ is the density of Gaussian random vector with the unit covariance matrix\footnote{of course, we could replace the unit covariance with any other positive definite covariance matrix, common for all distributions from the family in question} and expectation $\mu$. The space $\F$ is comprised of all affine functions on $\Omega=\bbr^m$, and because
\[
\ln\left(\int_{\bbr^m}e^{a^T\omega+b}p_\mu(\omega)d\omega\right)={b+a^T\mu+\half{a^T a}},
\]
Gaussian o.s. is good.
\paragraph{Poisson o.s.} Here $\Omega=\bbz_+^m$ is the discrete set of $m$-dimensional vectors with nonnegative integer entries, $P$ is the counting measure on $\Omega$,  $\M=\bbr^m_{++}$ is the set of positive $m$-dimensional vectors, and $p_\mu(\cdot)$, $\mu=[\mu_1;...;\mu_m]>0$, is the probability distribution of a random vector $\omega=[\omega_1;...;\omega_m]$, where $\omega_i$ is Poisson random variable with parameter $\mu_i$, and with independent of each other $\omega_1,...,\omega_m$. The space $\F$ is comprised of all affine functions on $\Omega$.
Note that
\[
\ln\left(\sum\limits_{\omega\in{\bbz}^m_+}
         \exp({a^T\omega+b})p_\mu(\omega)\right)={\sum_{i=1}^m(e^{a_i}-1)\mu_i}+b
\]
         is concave in $\mu$, and thus Poisson o.s. is good.
\paragraph{Discrete o.s.} Here $\Omega=\{1,...,m\}$ is a finite set, $P$ is a counting measure on $\Omega$, and $\M$ is comprised of all non-vanishing probability densities taken w.r.t. $P$, that is, the parameter space $\M$ is comprised of all $m$-dimensional vectors $\mu=[\mu_1;...;\mu_m]>0$ with entries summing up to 1, and the random variable $\omega$ distributed according to $p_\mu(\cdot)$ takes value $\omega\in\Omega$ with probability $\mu_\omega$. The space $\F$ is comprised of all real-valued functions on $\Omega=\{1,...,m\}$. Since
for $\phi\in \bbr^m$,
\[
\ln\left(\sum_{\omega\in \Omega}e^{\phi(\omega)} p_\mu(\omega)\right)= \ln\left(\sum_{\omega=1}^m e^{\phi_\omega} \mu_\omega\right)
\]
is concave in $\mu\in \M$, the Discrete o.s. is good.
\par
 We have seen that Gaussian, Poisson, and Discrete o.s.'s are good. More examples of good o.s.'s can be obtained by taking {\sl direct products}.
\paragraph{Direct products of good o.s.'s.} Let $\O_t=((\Omega_t,P_t),\{p_{\mu_t,t}(\cdot):\mu_t\in\M_t\},\F_t)$, $1\leq t\leq K$, be good o.s.'s. We can associate with this collection a new o.s. -- their {\sl direct product}, which, informally, describes the situation where our observation is $\omega^K=(\omega_1,...,\omega_K)$, with $\omega_t$ drawn, independently across $t$, from the o.s.'s $\O_t$. Formally, the direct product is the o.s.
{\small $$
\begin{array}{rcl}
\O^K&=&\big\{(\Omega^K,P^K)=(\prod\limits_{t=1}^K\Omega_t,P^K=P_1\times..\times P_K),\\
&&\{p_{\mu^K}(\omega_1,...,\omega_K)=p_{\mu_1,1}(\omega_1)...p_{\mu_K,K}(\omega_K):\mu^K=[\mu_1;...;\mu_K]\in\M^K=\M_1\times...\times \M_K\},\\
&&\multicolumn{1}{r}{\F^K=\{f(\omega_1,...,\omega_K)=\sum_{t=1}^Kf_t(\omega_t):f_t\in\F_t,1\leq t\leq K\}\big\},}\\
\end{array}
$$}
and this o.s. turns to be good provided all factor $\O_t$ are so. Note that with this definition, the parameter $\mu^K$ underlying the distribution of  observation $\omega^K$ is the collection of parameters $\mu_t$, $t=1,...,K$, underlying distributions of the components $\omega_1$,...,$\omega_K$ of $\omega^K$.
\par
Now consider the special case of this construction where all factors $\O_t$ are identical to each other:
$$
\O_t=((\Omega,P),\{p_{\mu}(\cdot):\mu\in\M\},\F),\,\,1\leq t\leq K.
$$
In this case  we can ``shrink'' the direct product $\O^K$ of our $K$ identical o.s.'s to the o.s. $\O^{(K)}$ by passing to observations $\omega^K=(\omega_1,...,\omega_K)$ with $\omega_1,...,\omega_K$ drawn, independently of each other,  from {\sl the same} density $p_{\mu}(\cdot)$, rather than being drawn from their own densities $p_{\mu_t,t}(\cdot)$. The formal description of $\O^{(K)}$ is as follows:
\begin{itemize}
\item for $\O^{(K)}$, the observation space $\Omega^K=\Omega\times...\times\Omega$ and the reference measure $P^K=P\times...\times P$ are exactly the same as for $\O^K$;
\item the family of probability densities  $\{p(\cdot)\}$ for $\O^{(K)}$ is $\{p^{(K)}_\mu(\omega_1,...,\omega_K)=\prod\limits_{t=1}^Kp_\mu(\omega_t):\mu\in\M\}$;
\item the family $\F$ for $\O^{(K)}$ is $\F^{(K)}=\{f^{(K)}(\omega_1,...,\omega_K)=f(\omega_1)+...+f(\omega_K):f\in\F\}$.
\end{itemize}
Note that $\O^{(K)}$ is a good o.s., provided $\O$ is so.
We shall refer to the just defined $\O^{(K)}$ as to {\sl stationary $K$-repeated observations} associated with $\O$.
\subsection{Pairwise hypothesis testing}\label{sect:pairwise}
\paragraph{Detectors and their risks.}
Let $\Omega$ be a Polish space and $\X_1,\X_2$ be two nonempty sets of Borel probability distributions on $\Omega$. Given a {\sl detector} -- a real-valued Borel function $\phi$ on $\Omega$, we define the {\sl risk} $\epsilon(\phi|\X_1,\X_2)$ of this detector w.r.t. $(\X_1,\X_2)$ as the smallest $\epsilon$ such that
\begin{equation}\label{eq1}
\begin{array}{rcl}
\int_\Omega \exp\{-\phi(\omega)\}p(d\omega)&\leq&\epsilon\,\,\forall p\in\X_1,\\
\int_\Omega \exp\{\phi(\omega)\}p(d\omega)&\leq&\epsilon\,\,\forall p\in\X_2,\\
\end{array}
\end{equation}
The sets $\X_1,\X_2$ give rise to two hypotheses, $H_1,\,H_2$, on the distribution of a random observation $\omega\in\Omega$, with $H_\chi$ stating that the distribution $p(\cdot)$ of this observation belongs to $\X_\chi$, $\chi=1,2$, while a detector $\phi$ gives rise to the test $T_\phi$ which, given an observation $\omega$, accepts $H_1$ (and rejects $H_2$) when $\phi(\omega)\geq0$, and rejects $H_1$ (and accepts $H_2$) otherwise. We define {\em the risk of a test} $T$ deciding between $H_1$ and $H_2$ as a maximal $p(\cdot)$-probability of the test rejecting the hypothesis $H_\chi$, $\chi=1,2$, when it is true:
$$\hbox{Risk}=\max\left[\sup_{p(\cdot)\in\X_1}p(\{\omega: T(\omega)=-1\}),\sup_{p(\cdot)\in\X_2}p(\{\omega: T(\omega)=1\})\right],$$
where $T(\omega)=1$ when the test $T$, as applied to observation $\omega$, accepts $H_1$, and $T(\omega)=-1$ otherwise.
Clearly, the risk of the test $T_\phi$ associated with detector $\phi$ is $\leq \epsilon:=\epsilon(\phi|\X_1,\X_2)$.
 Indeed, denoting by $p(\cdot)$ the distribution of the observation, when $H_1$ is true, the test rejects this hypothesis  when $\phi(\omega)<0$, and $p(\cdot)$-probability of this event, by the first inequality in (\ref{eq1}) and due to the definition of $\epsilon:=\epsilon(\phi|\X_1,\X_2)$, is at most $\epsilon$; when $H_2$ is true, the test rejects $H_2$ and accepts $H_1$ when $\phi(\omega)\geq0$, and $p(\cdot)$-probability of the latter event is  $\leq\epsilon$ due to the second inequality in (\ref{eq1}).

Note that if we have at our disposal a test, say $\overline{T}$, which decides between $H_1$ and $H_2$ with the risk bounded with $ \bar{\epsilon}\in(0,1/2)$, we can associate with it the detector
\[
\bar{\phi}(\omega)=\half\ln\left({1-\bar{\epsilon}\over\bar{\epsilon}}\right)\overline{T}(\omega)
\]
One can easily see that the risk of $\bar{\phi}(\cdot)$ satisfies the bounds of \rf{eq1}
with
\[
\epsilon=2\sqrt{\bar{\epsilon}(1-\bar{\epsilon})}.
 \]
As we will see in an instant, tests associated with detectors satisfying \rf{eq1} allow for a simple calculus -- one can ``propagate'' the tests properties to the case of repeated observations and multiple testing.
Our first goal is to describe a systematic construction of detectors satisfying \rf{eq1} in the situation where the underlying o.s. is good.

\paragraph{Near-optimal tests.} As far as testing pairs of hypotheses is concerned, the main result of \cite{GJuN2014} -- and the starting point of our developments  in this paper -- is as follows:
\begin{theorem}\label{the1} {\rm \cite[Theorem 2.1]{GJuN2014}} Let $\O=((\Omega,P),\{p_\mu(\cdot):\mu\in\M\},\F)$  be a good o.s., and let $X_1,X_2$ be two nonempty convex compact subsets of $\M$. Consider the optimization problem
\begin{equation}\label{eq2}
\Opt=\Opt(X_1,X_2,\O):=\max_{\mu\in X_1,\nu\in X_2}\left\{\psi(\mu,\nu):=\ln\left(\int_\Omega \sqrt{p_\mu(\omega)p_\nu(\omega)}P(d\omega)\right)\right\}
\end{equation}
 The function $\psi(\mu,\nu):\M\times\M\to\bbr$ is concave and continuous, so that {\rm (\ref{eq2})} is a convex optimization problem. This problem is solvable, and every optimal solution $(\mu_*,\nu_*)$ to this problem gives rise to the detector
 \be
 \phi(\omega)={1\over 2}\ln(p_{\mu_*}(\omega)/p_{\nu_*}(\omega)):\Omega\to\bbr
 \ee{dephi1}
 with the following properties:
\par
{\rm (i) [risk bound]} Denoting by $\X_\chi$ the set of all probability distributions on $\Omega$ with densities, w.r.t. $P$, of the form $p_\mu(\cdot)$ with $\mu\in X_\chi$, $\chi=1,2$, the risk of the detector $\phi$ w.r.t. $(\X_1,\X_2)$ is
$$
\epsilon(\phi|\O,\X_1,\X_2):=\epsilon_\star(\O,X_1,X_2)=\exp\{\Opt\};
$$
consequently, the risk of the test $T_\phi$ induced by the detector $\phi$, when deciding on the hypotheses $H_1$, $H_2$ associated with $\X_1$, $\X_2$, as at most $\exp\{\Opt\}$.\par
{\rm (ii) [near-optimality]} Assume that for some $\epsilon\in(0,1/4)$, ``in the nature'' there exists a test $T$ (a deterministic Borel function $T(\omega)$ of $\omega\in\Omega$ taking values 1 (``$H_1$ accepted, $H_2$ rejected'') and -1 (``$H_2$ accepted, $H_1$ rejected'') with risk in deciding on $H_1$, $H_2$ based on an observation $\omega$ at most $\epsilon$. Then
$$
\epsilon_*(X_1,X_2)\leq 2\sqrt{\epsilon}.
$$
\end{theorem}
We interpret (ii) as a statement of near-optimality of the test $T$ defined in (i) -- whenever the hypotheses $H_1$, $H_2$ associated with
$\X_1$, $\X_2$ can be decided upon with small risk, the risk of our test $T_\phi$ also is small. Note that (ii) also suggests that (maximal) degradation of the test performance when passing from the optimal test to the near-optimal one associated with detector \rf{dephi1} may be significant in the setting of Theorem \ref{the1}, when the decision is taken on the basis of one observation. When repeated observations are available, the sub-optimality of the proposed tests is expressed by a moderate absolute factor, when measured in terms of the length of the observation sample necessary to attain the desired testing accuracy.
\begin{corollary}\label{cor1} {\rm \cite[Proposition 2.1]{GJuN2014}} Let $\O=((\Omega,P),\{p_\mu(\cdot):\mu\in\M\},\F)$, $X_1,X_2$ be as in Theorem \ref{the1}, and $K$ be a positive integer. Assume we have at our disposal stationary $K$-repeated observation $\omega^K$ stemming from $\O$, and let $H_\chi$, $\chi=1,2$, be the hypotheses on the density  $p_\mu(\omega_1,...,\omega_K)=\prod_{t=1}^Kp_\mu(\omega_t)$ of $\omega^K$ stating that $\mu\in X_\chi$. Then
\item{\rm (i)} The optimal solution $(\mu_*,\nu_*)$ to the problem {\rm (\ref{eq2})} associated with $\O$, $X_1$, $X_2$, is
optimal for the same problem associated with $\O^{(K)},X_1,X_2$, and
$$
\epsilon_\star(\O^{(K)},X_1,X_2)=\left[\epsilon_\star(\O,X_1,X_2)\right]^K;
$$
moreover, the detectors $\phi$ and $\phi^{(K)}$ associated with this optimal solution by Theorem \ref{the1} as applied to $(\O,X_1,X_2)$ and $(\O^K,X_1,X_2)$, respectively, are linked by the relation
$$
\phi^{(K)}(\omega_1,...,\omega_K)=\sum_{t=1}^K\phi(\omega_t).
$$
\item{\rm (ii)}  Assume that for some $\epsilon\in(0,1/4)$ and some positive integer $\bar{K}$, there exists a test which decides on the hypotheses $H_1$, $H_2$ via stationary $\bar{K}$-repeated observation $\omega^{\bar{K}}$ with risk $\leq\epsilon$. Setting
$$
K=\hbox{\rm Ceil}\left({2\bar{K}\over 1-{2\ln(2)\over\ln(1/\epsilon)}}\right),
$$
we ensure that the risk of the test based on detector $\phi^{(K)}(\cdot)$ when deciding on the hypotheses $H_\chi$, $\chi=1,2$, is $\leq\epsilon$ as well. Note that $K/\bar{K}\to 2$ as $\epsilon\to+0$.
\end{corollary}
\subsection{Testing multiple hypotheses}\label{sect:multiple}
The completing ``building block'' from \cite{GJuN2014} we need is a simple technique for passing from pairwise tests to tests deciding on $N\geq 2$ hypotheses.
\paragraph{The situation} we are interested in is as follows. We are given a Polish observation space $\Omega$ along $N\geq2$ nonempty families $\X_j$, $1\leq j\leq N$, of Borel probability distributions on $\Omega$, and have at our disposal {\sl pairwise detectors} -- real-valued Borel functions $\phi_{ij}(\omega):\Omega\to\bbr$ and {\sl risk bounds} $\epsilon_{ij}\in(0,1]$, $1\leq i,j\leq N$,  such that
\begin{equation}\label{eq10}
\begin{array}{l}
\phi_{ij}(\cdot)=-\phi_{ji}(\cdot),\,\epsilon_{ji}=\epsilon_{ij},\,1\leq i,j\leq N\\
\int_\Omega \exp\{-\phi_{ij}(\omega)\}p(d\omega)\leq\epsilon_{ij}\,\,\forall p(\cdot)\in\X_i,\,1\leq i,j\leq N.
\end{array}
\end{equation}
Our goal is, given a positive integer $K$, to decide via stationary $K$-repeated observation $\omega^K=(\omega_1,...,\omega_K)$, with $\omega_1,...,\omega_K$ drawn, independently of each other, from a distribution $p(\cdot)\in\bigcup\limits_{j=1}^M\X_j$, between hypotheses $H_j$, $1\leq j\leq N$, with $H_j$ stating that $p(\cdot)\in\X_j$.
\par
Note that the present setting is a straightforward extension of the situation considered in section \ref{sect:pairwise}.
In particular, given what was called in this section ``a detector $\phi$ with $(X_1,X_2)$-risk $\leq \epsilon$,'' and setting $\phi_{1,2}=\phi$, $\phi_{2,1}=-\phi$, $\phi_{1,1}\equiv\phi_{2,2}\equiv 0$, $\epsilon_{1,2}=\epsilon_{2,1}=\epsilon$, $\epsilon_{1,1}=\epsilon_{2,2}=1$, we meet the requirements (\ref{eq10}) corresponding to the case $N=2$. Vice versa, given detectors and risks satisfying (\ref{eq10}) for the case $N=2$ and setting $\phi=\phi_{1,2}$, $\epsilon=\epsilon_{1,2}$, we get a detector $\phi$ with $(X_1,X_2)$-risk $\leq\epsilon$.
\par
We are about to ``aggregate'' detectors $\phi_{ij}$ into a testing procedure deciding on the hypotheses $H_1,...,H_N$ via stationary $K$-repeated observations $\omega^K$. It makes sense to consider a slightly more general problem, specifically, as follows: assume that on the top of the setup data $\X_j$, $\phi_{ij}(\cdot)$, $\epsilon_{ij}$, we are given an $N\times N$ symmetric (proximity) matrix $\C$ with zero-one entries and zero diagonal. We interpret the relation $\C_{ij}=0$ as {\sl closeness} of hypotheses $H_i$ and $H_j$\footnote{In the general case considered in \cite{GJuN2014}, the closeness (i.e., the zero-one matrix $\C$ with zero diagonal) not necessary is symmetric; here we restrict ourselves with symmetric closeness only.}, and we refer to indices $i$, $j$ (and hypotheses $H_i$, $H_j$) such that $\C_{ij}=0$ as {\sl $\C$-close.} We expect our testing procedure not to reject the true hypothesis, while rejecting all hypotheses which are not {\sl not} $\C$-close to it. On the other hand, we do not care about distinguishing the true hypothesis from  $\C$-close alternatives.
\paragraph{The construction.} Given $\X_j$, $\phi_{ij}(\cdot)$, $\epsilon_{ij}$ satisfying (\ref{eq10}) along with positive integer $K$, let us set
$$
\phi_{ij}^K(\omega^K)=\sum_{t=1}^K\phi_{ij}(\omega_t),\,\,\epsilon_{ij}^{(K)}=\epsilon_{ij}^K,
$$
and let $\X_j^K$, $j=1,...,N$ be a family of probability distributions of $\omega^K=(\omega_1,...,\omega_K)$ where $\omega_t$ are i.i.d. with distribution $p(\cdot)\in \X_j$. Clearly, (\ref{eq10}) implies that
\[
\begin{array}{l}
\phi_{ij}^K(\cdot)=-\phi_{ji}^K(\cdot),\,\epsilon_{ji}^{(K)}=\epsilon_{ij}^{(K)},\,1\leq i,j\leq N\\
\int_{\Omega^K} \exp\{-\phi_{ij}^K(\omega^K)\}p(d\omega^K)\leq\epsilon_{ij}^{(K)}\,\,\forall p(\cdot)\in\X_i^K,\,1\leq i,j\leq N.
\end{array}
\]
Now, let $\alpha=[\alpha_{ij}]$ be a skew-symmetric $N\times N$ matrix, and let
$$\overline{\phi}_{ij}(\omega^K)=\phi_{ij}^K(\omega^K)-\alpha_{ij}.$$
We associate with $\overline{\phi}$ the test $\T_K$ which, given observation $\omega^K$, builds the $N\times N$ matrix with the entries $\overline{\phi}_{ij}(\omega^K)$ and accepts all hypotheses $H_i$ which satisfy the condition
$$
\overline{\phi}_{ij}(\omega^K)>0\;\;\forall (j: C_{ij}=1),
$$
and rejects all remaining hypotheses. Note that $\T_K$ can accept no hypotheses at all, or can accept more than one hypothesis.
\par
The properties of $\T_K$ are summarized in the following simple statement (see \cite[section 2.3.1]{GJuN2014}):
\begin{proposition}\label{propmult}Let
 \begin{equation}\label{eq32}
 \varepsilon=\max\limits_{1\leq i\leq N} \sum_{j:\C_{ij}=1}\epsilon_{ij}^K\exp\{-\alpha_{ij}\}.
 \end{equation}
Let $\omega^K=(\omega_1,...,\omega_K)$ be sampled independently from the distribution $p_*(\cdot)\in \X_{i_*}$, for some ${i_*}\in \{1,...,N\}$.      Then
{the $p_*(\cdot)$-probability of the event ``$\T_K$ does not accept the true hypothesis $H_{i_*}$ or accepts a hypothesis $H_i$ which is not $\C$-close to $H_{i_*}$'' is $\leq\varepsilon.$}
\end{proposition}
The risk bound $\varepsilon$, as given by (\ref{eq32}), depends on the ``shifts'' $\alpha_{ij}=-\alpha_{ji}$, and we would like to use the shifts resulting in as small as possible value of $\varepsilon$ in \rf{propmult}. It is shown in \cite[section 3]{GJuN2014} that the corresponding shifts solve
the convex optimization problem
\begin{equation}\label{eq33}
\Opt=\min_{\alpha\in\bbr^{N\times N}}\left\{f(\alpha)=\max\limits_{1\leq i\leq N} \sum_{j:\C_{ij}=1}\epsilon_{ij}^K\exp\{\alpha_{ij}\}:\alpha=-\alpha^T\right\},
\end{equation}
the optimal $\varepsilon$ in \rf{eq32} being exactly $\Opt$. Moreover, from \cite[Proposition 3.3]{GJuN2014} it immediately follows (see below), that $\Opt$ is nothing but the spectral norm $\|D\|_{2,2}$ of the entry-wise nonnegative symmetric matrix
$$
D=[d_{ij}=\epsilon_{ij}^K\C_{ij}]_{1\leq i,j\leq N}.
$$
Moreover, it is immediately seen (cf. \cite[section 3.2]{GJuN2014}) that {\sl if the Perron-Frobenius eigenvector $g$ of the entrywise nonnegative symmetric matrix $D$ is positive, an optimal solution to (\ref{eq33}) is given by}
\begin{equation}\label{eq34}
\alpha_{ij}=\ln(g_j)-\ln(g_i).
\end{equation}
In the general case the Perron-Frobenius eigenvector of $D$ can have zero entries, and, moreover, (\ref{eq33}) may happen not to have optimal solutions at all; we, however, can easily find ${\delta}$-optimal, with a whatever small ${\delta}>0$, solutions to the problem by utilizing in (\ref{eq34}) in the role of $g$ the Perron-Frobenius eigenvector of a matrix $D'$ with the entries $d'_{ij}>d_{ij}$ close to $d_{ij}$, specifically, such that $\|D'\|_{2,2}\leq \|D\|_{2,2}+{\delta}$.
\begin{quote}
{\small In fact, Proposition 3.3 from \cite{GJuN2014} states that if $E=[e_{ij}]$ is a symmetric $N\times N$ matrix with zero diagonal and {\em positive} off-diagonal entries, then the optimal value in the optimization problem
$$
\Opt_E=\min_\alpha\left\{\max_{1\leq i\leq N} \sum_{j=1}^Ne_{ij}\exp\{\alpha_{ij}\}:\alpha=-\alpha^T\right\}
$$
is equal to the spectral norm $\|E\|_{2,2}$.
Let now $e_{ij}=d_{ij}(=\epsilon_{ij}^K\C_{ij})$ for $1\leq i,j\leq N$ such that $\C_{ij}=1$, and let $e_{ij}={\delta}>0$ for $i\neq j$ such that $\C_{ij}=0$, so that all off-diagonal entries of $E$ are positive. Note that when ${\delta} \downarrow 0$, both the spectral norms of $\|E\|_{2,2}$ and $\|D\|_{2,2}$, and  the optimal values $\Opt$ and $\Opt_E$ become arbitrarily close to each other. Since for the ``perturbed'' matrix $E$ the spectral norm $\|E\|_{2,2}$ coincides with the optimal value $\Opt_E$, the same holds true for the ``unperturbed'' matrix $D$.}
\end{quote}
\section{Sequential Hypothesis Testing}\label{sec:seq1}
\subsection{Problem setting}
Let us consider the situation of section \ref{sect:multiple}. Specifically, assume that we are given $N\geq 2$ nonempty families $\X_j$ of Borel probability distributions on a Polish observation space $\Omega$.
Let, further,
   $$
   \J:=\{1,2,...,N\}=\bigcup\limits_{i=1}^I\J_i
   $$
   be a partition of the set of indices of $X_j$'s into $I\geq 2$ non-overlapping nonempty groups $\J_1,...,\J_I$.
We associate with the sets $\X^i=\bigcup_{j\in \J_i} \X_j$ the hypotheses $H_i$, stating that the elements $\omega_1,...,\omega_K$ of the $K$-repeated stationary observation sample $\omega^K$ are drawn independently from the common distribution $p\in \X^i$, and our goal is to decide from observation $\omega^K$ on the hypotheses $H_1,...,H_I$.
\par
 It is convenient to think about the (values of the) indices $i=1,...,I$ of the sets $\X^i$ as of the {\sl colors} of these sets.\footnote{Let us agree that these colors are inherited by different entities associated with $\X^i$'s such as index groups $\J_i$, indices $j\in \J_i$, corresponding sets $\X_j$ and distributions $p\in \X^i$.}
 \par
 We also assume that we have at our disposal pairwise detectors -- Borel functions $\phi_{ij}: \Omega\to \bbr$ -- for the sets $\X_j$, $j=1,...,N$ which satisfy relations \rf{eq10}.
From now on, we make the following assumption:
\begin{quote}
   {\bf A:} {\sl When the indices $j$ and $j'\in\J$ are of different colors (in other words, do not belong to the same group $\J_i$), the risk $\epsilon_{jj'}(=\epsilon_{j'j})$ of the detector $\phi_{jj'}$ (same as the detector $\phi_{j'j}$) satisfies $\epsilon_{jj'}<1$.
   }
   \end{quote}
   Note that assumption {\bf A} implies that for $j$ and $j'$ of different colors
the sets of distributions
$\X_j$ and $\X_{j'}$ are at positive Hellinger distance from each other.
This implies that our color assignments are unambiguous, and our goal may be reformulated as that of identifying the color of the distribution underlying the observations.
\par
 From the results of section \ref{sect:multiple} it easily follows that if assumption {\bf A} holds then for any given  $\epsilon>0$, we can decide on the hypotheses $H_1,...,H_I$ with risk $
\leq\epsilon$ (meaning that the probability to reject the true hypothesis, same as the probability to accept a wrong one, is $\leq\epsilon$), provided that the number $K$ of observations is large enough.
This being said, the ``large enough'' $K$ could be indeed quite large { if there are pairs $j\in \J_i,j'\in\J_{i'}$, $i\neq i'$, with close to 1  values of $\epsilon_{jj'}$}. As a tradeoff, we can switch from decision rules based on $K$ observations to {\sl sequential} decision rules, for which the decision is made on the basis of on-line adjustable number of observations. We can expect that if we are lucky and the distribution $p_*$ underlying our observation is ``deeply inside'' of some $\X^{i_*}$ and thus is ``far'' from all $\X^i$, $i\neq i_*$, the true hypothesis $H_{i_*}$ will be accepted much sooner than in the case when $p_*$ is close to some of ``wrong'' $\X^i$'s. Our objective now is  to build sequential tests utilizing the results of section \ref{sec:oldies}.

\subsection{Sequential test: construction}\label{sect:construction}
\paragraph{The setup} for our ``generic'' sequential test is given by
\begin{enumerate}
\item required {\sl risk} $\epsilon\in(0,1)$;
\item positive integer $S$ -- {\sl number of stages}, along with the following entities, defined for $1\leq s\leq S$ and forming {\sl $s$-th component} of the setup:
\begin{enumerate}
\item positive reals $\epsilon_s$, $1\leq s\leq S$, such that ${ \sum_{s=1}^S\epsilon_s}=\epsilon$;
\item representations $\X_j=\bigcup\limits_{\iota=1}^{\iota_{js}} \X_{j\iota s}$, $j\in \J$, where $\X_{j\iota s}$ are nonempty  subsets of $\X_j$;\\
\item {\sl tolerances} $\delta_s\in(0,1)$.
\end{enumerate}
\end{enumerate}
To avoid messy notation, we enumerate, for every $s\leq S$, the sets $\X_{j\iota s}$, $1\leq j\leq N, 1\leq \iota\leq\iota_{js}$, and call the resulting sets $\Z_{1s},...,\Z_{L_ss}$. Thus, $\Z_{qs}$ is one of the sets $\X_{j\iota s}$, and we assign  to $\Z_{qs}$ and the index $q$ the same color as that of $j$.
\paragraph{Detectors.}
We suppose for every $s\leq S$ and every pair $(q,q')$, $1\leq { q,q'}\leq L_s$, we are given pairwise detectors $\phi_{qq',s}$ and reals $\epsilon_{qq',s}\in(0,1]$ associated with $\Z_{qs}$ and $\Z_{q's}$ and satisfying the relations
\[\begin{array}{l}
\phi_{qq',s}(\cdot)=-\phi_{q'q,s}(\cdot),\,\epsilon_{qq',s}=\epsilon_{q'q,s},\,1\leq q,q'\leq L_s\\
\int_\Omega \exp\{-\phi_{qq',s}(\omega)\}p(d\omega)\leq\epsilon_{qq',s}\,\,\forall p(\cdot)\in\Z_{qs},\,1\leq q,q'\leq L_s.
\end{array}
\]

\paragraph{$s$-closeness.}
Let $H_{qs}$ be the hypotheses on the distribution $p(\cdot)$ of an observation stating $p(\cdot)\in \Z_{qs}$; by convention, $H_{qs}$ is of the same color as $\Z_{qs}$. Let us say that hypothesis $H_{q's}$ is {\sl $s$-close} to hypothesis $H_{qs}$ (same as $q'$ is $s$-close to $q$), if either $H_{qs}$ and $H_{q's}$ are of the same color, or $\epsilon_{qq',s}>\delta_s$ if  they are of different colors. Let  $\C=\C^s$ be the $L_s\times L_s$ matrix with $(q,q')$-entry equal to 0 if and only if $q$ is $s$-close to $q'$, and equal to 1 otherwise. This matrix clearly meets the requirements imposed on $\C$ in section \ref{sect:multiple}:
it is a symmetric matrix (recall that $\epsilon_{qq',s}=\epsilon_{q'q,s}$) with 0/1 entries and zero diagonal.

\paragraph{Tests $T_s$.}
Now let us apply to the collection of hypotheses $\{H_{qs}:1\leq q\leq L_s\}$, detectors $\{\phi_{qq',s}(\cdot):1\leq q,q'\leq L_s\}$ and the just defined matrix $\C^s$ the construction from section \ref{sect:multiple}, assuming that when deciding upon hypotheses $H_{1s},...,H_{L_ss}$, we have at our disposal $k$-repeated observation $\omega^k$, with a given $k$. Specifically, consider the optimization problem
\begin{equation}\label{seq9}
\Opt(k,s)=\min_{\alpha\in\bbr^{L_s\times L_s}}\left\{f_s(\alpha):=\max_q\sum_{q':\,\C^s_{qq'}=1}\epsilon_{qq',s}^k\exp\{\alpha_{qq'}\}:\alpha=-\alpha^T\right\};
\end{equation}
as was shown in section \ref{sect:multiple}, the value $\Opt(k,s)$ is nothing but the spectral norm of the entry-wise nonnegative symmetric matrix
with the entries $\epsilon_{qq',s}^k\C^{s}_{qq'}$.
\par
Since $\epsilon_{qq',s}\leq \delta_s\in(0,1)$ when $\C^s_{qq'}=1$,  $\Opt(k,s)$ goes to 0 as $k\to\infty$, so that the smallest $k=k(s)$ such that $\Opt(k,s)<\epsilon_s$ is well defined. And since $\Opt(k(s),s)<\epsilon_s$, problem (\ref{seq9}) with $k=k(s)$, whether solvable or not, admits a feasible solution $\bar{\alpha}^{(s)}$ such that $f_s(\bar{\alpha}^{(s)})\leq\epsilon_s$. Applying to detectors $\overline{\phi}_{qq',s}(\cdot)=\phi_{qq',s}(\cdot)-\bar{\alpha}^{(s)}_{qq'}$ and to $\C=\C^s$ the construction from section \ref{sect:multiple}, we get a test $T_{s}$ deciding on the hypotheses $H_{qs}$, $1\leq q\leq L_s$, via $k(s)$-repeated observation $\omega^{k(s)}$, with properties as follows:
\begin{quote}{\em
Let $\omega_1,...,\omega_{k(s)}$ be drawn, independently of each other, from common distribution $p_*(\cdot)$ obeying the hypothesis $H_{q_*s}$ for some $1\leq q_*\leq L_s$. Then
 {the $p_*(\cdot)$-probability of the event ``$T_s$ does not accept the true hypothesis $H_{q_*s}$ or accepts a hypothesis $H_{qs}$ with $q$ not $s$-close to $i_*$'' is at most $\epsilon_s$.}
}
\end{quote}
\paragraph{$S$-th setup component.}
We assume that {\sl when $s=S$, the partition of $\X_j$ is trivial: $\iota_{jS}=1$ and ${\X}_{j1S}={\X}_j$ for all $j\leq N$.} Furthermore, we define $K=k(S)$ such that
 {\sl $\delta_S\geq\epsilon_{qq',S}$ whenever $q$, $q'$ are of different colors}, implying that {\sl $q$ and $q'$ are $S$-close if and only if $q$ and $q'$ are of the same color}.\footnote{We can  do so because by assumption {\bf A}   {\sl all quantities $\epsilon_{qq',s}$, $1\leq q,q'\leq L_S=N$, with $q$, $q'$ of different colors are less than 1.}}
For the reasons which will become clear in a moment, we are not interested in those components of our setup for which $k(s)>K$; if components with this property were present in our original setup, we can just eliminate them, reducing $S$ accordingly. Finally, we can reorder the components of our setup to make $k(s)$ nondecreasing in $s$. Thus, from now on we assume that
$$
k(1)\leq k(2)\leq...\leq k(S)=:K.
$$

\paragraph{Sequential test $\T$} corresponding to the outlined setup when applied to the observation $\omega^K$ works by stages $s=1,2,...,S$: at stage $s$ we apply test $T_s$ to the initial fragment $\omega^{k(s)}$ of $\omega^K$.
If the outcome of the latter test is acceptance of a nonempty set of hypotheses $H_{qs}$ {\sl and all these hypotheses are of the same color $i$}, test $\T$ accepts the hypothesis $\H_i$ and terminates, otherwise it proceeds to stage $s+1$ (when $s<S$) or terminates without accepting any hypothesis $(s=S)$.\par
Note that by construction $\T$ never accepts more than one of the hypotheses $H_1,...,H_I$.
\subsection{Sequential test: analysis}\label{sectriskanal} For a distribution $p\in \X=\bigcup\limits_{j=1}^N \X_j$, let $s[p]\in\{1,2,...,S\}$ be defined as follows: for every $s$, $p$ belongs to (perhaps, several) of the sets $\Z_{qs}$, $1\leq q \leq L_s$. Further, let $Q_s[p]$ be the set of all $q$'s such that $p\in \Z_{qs}$. Note that the color of every  $q\in Q_s[p]$ (recall that the sets $\Z_{qs}$ and the corresponding values of $q$  have already been assigned colors) is the same as the color of $p$. Now, given $p\in \X$, for some $s\leq S$ it may happen that
\begin{equation}\label{mayhappen}
\exists q\in Q_s[p]: \hbox{\ all $q'$ $s$-close to $q$ are of the same color as $q$.}
\end{equation}
In particular, the latter condition definitely takes place when $s=S$, since, as we have already seen, $q$ and $q'$ are $S$-close if and only if $q$ and $q'$ are of the same color. Hence for $s=S$ the conclusion in (\ref{mayhappen}) is satisfied for all $q\in Q_S[p]$. Now let $s[p]$ be the smallest $s$ such that (\ref{mayhappen}) takes place, and let the corresponding $q$ be denoted by $q[p]$. Thus for all $p\in \X$,
\[
\begin{array}{l}
s[p]\in\{1,2,...,S\},\,\, q[p]\in\{1,...,L_{s[p]}\},\,\, p\in {\Z}_{q[p]s[p]};\\
\hbox{whenever $q'\in\{1,...,L_{s[p]}\}$ is $s[p]$-close to $q[p]$, $q'$ and $q[p]$ are of the same color.}\\
 \end{array}
\]
The main result of this section is as follows.
\begin{proposition}\label{prop_sequential} Let $\omega_1,...,\omega_K$ be drawn, independently of each other, from a distribution $p_*\in \X(=\bigcup_{j=1}^N \X_j)$, so that $p_*\in \X_{j_*}$ for some $j_*\leq N$, and let $i_*$ be such that $j_*\in \J_{i_*}$ (i.e., $i_*$ is the color of $p_*$). Then $p_*$-probability of the event
$$
\E=\left\{\omega^K: \hbox{\rm \begin{tabular}{l}$\T$, as applied to $\omega^K$, terminates not later than at the stage $s[p_*]$\\
 and accepts upon termintation the true hypothesis $\H_{i_*}$\\
 \end{tabular}} \right\}
$$
is at least $1-\epsilon$.
\end{proposition}

\section{Implementing sequential test}\label{sec:impl1}
Observe that in order for the just described sequential test to  recover the ``non-sequential'' test from section \ref{sect:multiple}, it suffices to utilize setup with $S=1$ (recall that our construction fully specifies the setup component with $s=S$). This being said, with our sequential test the running time (the number of observations used to make the inference) depends on the true distribution underlying the observations, and we can use several degrees of freedom in our setup in order to save on the number of observations when the true distribution is ``deeply inside'' the true hypothesis. We are about to illustrate some of the options available for the basic good o.s.'s presented in section \ref{sec:good_os}.
\subsection{Preliminaries} Assume that our o.s. is either Gaussian, or Poisson, or Discrete, see section \ref{sect:examples}. We denote by $n$ the dimension of the associated parameter vector $\mu$: $\mu=[\mu_1;...;\mu_n]$.
\par
Assume that we are given the risk level $\epsilon\in(0,1)$ and a collection of $J$ nonempty compact convex sets $X_j\subset \M$ painted in $I\geq 2$ colors (i.e., the set of indices $\{1,...,J\}$ is split into $I\geq2$ non-overlapping nonempty sets $\J_1,...,\J_I$, $i$ being the common color of all sets $X_j$, $j\in\J_i$).
Sets $X_j$ give rise to the sets $\X_j$ of probability distributions defined by the corresponding densities $p_\mu(\cdot)$, $\mu\in X_j$. Note that for all considered o.s.'s
 different values of the parameter $\mu\in \M$ correspond to different probability densities $p_\mu$.
We assume that the sets $X_j$, $X_{j'}$ of different colors do not intersect thus giving rise to non-intersecting sets of distributions $\X_j$ and $\X_{j'}$. Given this input, we intend to specify the setup for sequential test deciding on the associated hypotheses $H_i$, $1\leq i\leq I$.
\par
Let $\psi(\mu,\nu):\,\M\times\M\to\bbr$ be the rate function
 \[
\psi(\mu,\nu)=\ln\left(\int_{\Omega}\sqrt{p_\mu(\omega)p_\nu(\omega)}P(d\omega)\right),\]
 associated with the o.s. in question.
 An immediate computation shows that (cf. \cite[section 2.3]{GJuN2014})
\[
\psi(\mu,\nu)=\left\{\begin{array}{ll}-{1\over 8}\|\mu-\nu\|_2^2,\,\mu,\nu\in\M=\bbr^n,&\hbox{Gaussian case}\\
-{1\over 2}\sum_{\ell=1}^n(\sqrt{\mu_\ell}-\sqrt{\nu_\ell})^2,\,\mu,\nu\in\M=\{x\in\bbr^n:x>0\},&\hbox{Poisson case}\\
\ln\left(\sum_{\omega=1}^n\sqrt{\mu_\omega\nu_\omega}\right),\,\mu,\nu\in\M=\{x\in\bbr^n:x>0,\sum_{\omega=1}^nx_\omega=1\},&\hbox{Discrete case}\\
\end{array}\right.
\]
Theorem \ref{the1} states that for the considered o.s.'s, given two nonempty convex compact subsets $X,Y$ of $\M$  and setting
$$
\Psi_{XY}=\max_{\mu\in X,\nu\in Y} \psi(\mu,\nu),
$$
the quantity $\exp\{\Psi_{XY}\}$ is exactly the risk of the detector $\phi(\cdot)$ yielded by Theorem \ref{the1} as applied to $X,Y$ and the o.s. in question. Besides this, $\psi(\mu,\nu)$ clearly is a smooth concave and symmetric $(\psi(\mu,\nu)=\psi(\nu,\mu))$ function on its domain, and $\psi(\mu,\nu)<0$ whenever $\mu\neq\nu$.
\par
For $j\in \{1,..., J\}$, we set
$$
\psi_j(\mu)=\max_{\nu\in X_j} \psi(\mu,\nu):\M\to\bbr.
$$
Since $\psi(\mu,\nu)$ is concave in $\mu,\nu\in\M$, and $X_j$ is a convex compact set, the functions $\psi_j(\cdot)$ are concave and continuous on $\M$.
  \par
Let now $\cR$ be the set of all ordered pairs $(j,j')$, $1\leq j,j'\leq J$, {\sl with $j,j'$ of different colors;} note that $(j,j')\in\cR$ if and only if $(j',j)\in\cR$.
\par
 For a pair $j,j'$, $1\leq j,j'\leq J$, let
\begin{equation}\label{seq300}
\psi_{jj'}=\Psi_{X_jX_{j'}}:=\max_{\mu\in X_j,\nu\in X_{j'}}\psi(\mu,\nu)=\max_{\mu\in X_j}\psi_{j'}(\mu).
\end{equation}
                    Note that for $(j,j')\in\cR$ the convex compact sets $X_j$, $X_{j'}$ do not intersect, thus $\psi_{jj'} <0$, and  the objective in the optimization problem on the right hand side of \rf{seq300} is negative on the compact feasible set of the problem. We set
\be
d=\min_{(j,j')\in\cR}[-\psi_{jj'}],
\ee{defd}
so that $d>0$.
\par
Finally, for  $(j,j')\in \cR$ and a nonnegative $r$
we say that a linear inequality $\ell(\mu)\leq 0$, $\mu\in \M$, defines a {\sl $(jj'r)$-cut} if for all $\mu\in X_j$ such that $\ell(\mu)\leq 0$,
\[\psi_{jj'}(\mu)\leq -r.\]
{\sl Example: default cuts.} For $(j,j')\in\cR$, let $(\mu_{jj'},\nu_{jj'})$ be $(\mu,\nu)$-components of an optimal solution to the optimization problem (\ref{seq300}).
Setting
$$
e_{jj'} =\nabla_\mu\psi(\mu_{jj'},\nu_{jj'}), \,\,f_{jj'}=\nabla_\nu\psi(\mu_{jj'},\nu_{jj'}),\,\,(j,j')\in \cR,
$$
and invoking optimality conditions for (\ref{seq300}) along with concavity of $\psi(\cdot,\cdot)$, we get for all $(j,j')\in\cR$:
\begin{equation}\label{seq301}
\forall (\mu\in X_j,\nu\in X_{j'}):\;\left\{\begin{array}{ll}
e_{jj'}^T[\mu-\mu_{jj'}]\leq 0,&(a)\\
f_{jj'}^T[\nu-\nu_{jj'}]\leq 0,&(b)\\
\psi(\mu,\nu)\leq \psi_{jj'} +e_{jj'}^T[\mu-\mu_{jj'}]+f_{jj'}^T[\nu-\nu_{jj'}].&(c)\\
\end{array}\right.
\end{equation}
We conclude that setting
$$
\ell_{jj'}(\mu)=\psi(\mu_{jj'},\nu_{jj'})+e_{jj'}^T(\mu-\mu_{jj'})-r,
$$
we get an affine function of $\mu\in\M$ which upper-bounds $\psi_{jj'}(\cdot)-r$ on $X_j$.
\begin{quote} {\small Indeed,
for any $\nu\in X_{j'},\mu\in X_j$
\bse
\psi(\mu,\nu)&\leq& \psi({\mu}_{jj'},{\nu}_{jj'})+e_{jj'}^T(\mu-{\mu}_{jj'}) +f_{jj'}^T(\nu-{\nu}_{jj'})\\
&\leq& \psi({\mu}_{jj'},{\nu}_{jj'})+e_{jj'}^T(\mu-{\mu}_{jj'})
\ese
(we have used (\ref{seq301}.$c$,$b$)). Taking in the resulting inequality the supremum over $\nu\in X_{j'}$, we arrive at $\ell_{jj'}(\mu)\geq\psi_{j'}(\mu)-r$, $\mu\in X_j$.
}
\end{quote}
The bottom line is that $\ell_{jj'}(\mu)\leq0$ is a $(jj'r)$-cut; we shall refer to this cut as {\sl default}.
{
\par\noindent
{\sl Default cuts: refinement in Gaussian case.} As it will become clear from the sequel, for our purposes the  larger is the set
$
\{\mu\in X_j:\ell(\mu)\leq0\}$, the better is $(jj'r)$-cut $\ell(\mu)\leq0$. It is immediately seen that in the Gaussian case the default cuts as defined above can be easily improved. Indeed, let Gaussian case take place, and let $(j,j')\in\cR$. In this case, in the above notation we have $\psi(\mu,\nu)=-{1\over 8}\|\mu-\nu\|_2^2$, $(a_{jj'},b_{jj'})$ is a pair of $\|\cdot\|_2$-closest to each other points in $X_j$ and $X_{j'}$, $e_{jj'}=-{1\over 4}(a_{jj'}-b_{jj'})$ is nonzero and separates $X_j$ and $X_{j'}$:
$$
\forall (\mu\in X_j,\nu\in X_{j'}): e_{ij}^T(\mu-\nu)\leq e_{ij}^T(a_{jj'}-b_{jj}')=-{1\over 4}\delta_{jj'}^2,\,\,\delta_{jj'}=\|a_{jj'}-b_{jj'}\|_2,
$$
whence, setting $\eta_{jj'}=-e_{jj'}/\|e_{jj'}\|_2$,
\begin{equation}\label{whatisetacut}
\psi_{j'}(\mu)=-{1\over 8}\min_{\nu\in X_{j'}}\|\mu-\nu\|_2^2\leq -{1\over 8}\min_{\nu\in X_{j'}}|\eta_{jj'}^T[\mu-\nu]|^2\ \&\ (\forall \mu\in X_j): \eta_{jj'}^T(\mu-a_{jj'})\geq0.
\end{equation}
It follows that
$$
\mu\in X_j\Rightarrow \psi_{j'}(\mu)\leq -{1\over 8}\min_{\nu\in X_{j'}}|\eta_{jj'}^T[\mu-\nu]|^2=-{1\over 8}\big[\underbrace{\eta_{jj'}^T(\mu-a_{jj'})}_{\geq0}+\delta_{jj'}\big]^2,
$$
implying that the linear inequality
\begin{equation}\label{implyingcut}
\ell_{jj'r}(\mu):=2\sqrt{2r}-\delta_{jj'}-\eta_{jj'}^T(\mu-a_{jj'})\leq0
\end{equation} is a $(jj'r)$-cut.
On the other hand, the default $(jj'r)$ cat as defined by our general construction in Gaussian case is
$$
r-{1\over 8}\delta_{jj'}^2-{\delta_{jj'}\over 4}\eta_{jj'}^T(\mu-a_{jj'})\leq0.
$$
Taking into account that, as we have already mentioned, $X_j\subset \{\mu:\eta_{jj'}^T(\mu-a_{jj'})\geq0\}$, it is immediately seen that the validity of the latter inequality at $\mu\in X_j$ implies the validity of (\ref{implyingcut}), that is, the default cut cuts off $X_j$ a smaller set than the cut (\ref{implyingcut}). By this reason, {\sl from now on, when in the Gaussian case, we refer to the $(jj'r)$-cut {\rm (\ref{implyingcut})} as to the default one.}
}
\subsection{Specifying the setup}
The setup for our sequential test is as follows.
\begin{enumerate}
\item We select a sequence of positive integers $\{\bar{k}(s)\}_{s=1}^\infty$ satisfying
\begin{equation}\label{monotonicity}
\bar{k}(1)=1,\,\,\bar{k}(s)<\bar{k}(s+1)\leq 2\bar{k}(s),\,\,s=1,2,...
\end{equation}
and specify $S$ as the smallest positive integer such that
\begin{equation}\label{seq10}
\bar{k}(S)>d^{-1}\ln\left(SJ^2/\epsilon\right);
\end{equation}
($S$ is well defined due to $\bar{k}(s)\geq s$).
\par
For $1\leq s\leq S$, we set
\begin{equation}\label{seq11}
\begin{array}{c}
\epsilon_s={{\epsilon\over S}},\;
r(s)=\bar{k}(s)^{-1}\ln(SJ^2/\epsilon),\;\;\delta_s=\exp\{-r(s)\}.\\
\end{array}
\end{equation}
\item
For every $j\in \{1,...,J\}$ and every $s\in \{1,..., S\}$, we specify  closed convex subsets $X_{j\iota s}$, $1\leq \iota\leq \iota_{js}$, of $X_j$ as follows. For every pair $(j,j')\in\cR$, we select somehow
a $(jj'r(s))$-cut $\ell_{jj's}(\cdot)\leq0$ and set
\begin{equation}\label{partition}
\begin{array}{rcl}
X^{j'}_{js}&=&\{\mu\in X_j: \ell_{jj's}(\mu) \geq 0\},\; j'\in\overline{\J}_j,\\
X^j_{js}&=&\{\mu\in X_j: \ell_{jj's}(\mu)\leq0,\,\, j'\in \overline{\J}_j\},
\end{array}
\end{equation}
where for  $1\leq j\leq J$ the set $\overline{\J}_j$ contains all indices $1\leq j'\leq J$
of the color different from that of $j$. Eliminating form this
 list all sets which are empty,  we end up with a number $\iota_{js}\leq J$ of {\sl nonempty} convex compact sets $X_{j\iota s}$, $1\leq\iota\leq\iota_{js}$, with $X_j$ being their union.
 \par
Observe that $r(S)<d$ by (\ref{seq10}), and for ($j,j')\in\cR$ we clearly have
$$
\max_{\mu\in X_j}\psi_{j'}(\mu)=\psi_{jj'}\leq \max_{(j,j')\in \cR}\psi_{jj'}= -d
$$
implying that {affine functions} $\ell_{jj'S}(\mu)\equiv -1$ {specify legitimate $(jj'r(S))$-cuts}. These are exactly the cuts we use when $s=S$.
\end{enumerate}
Let us agree that for a subset of the parameter space $\M$ denoted by a capital Latin letter, the script version of the letter denotes the set of all densities $p_\mu$ with parameters  $\mu$ from the original set; e.g., if $A$ denotes a subset of $\M$, then $\A=\{p_\mu:\mu\in A\}$. This agreement works in both directions: if, say, $\A$ is a subset of $\{p_\mu:\mu\in M\}$, then $A=\{\mu\in\M:p_\mu\in\A\}$.
With this convention, the above sets $X^{j'}_{js}$, $j\in\overline{\J}_j$,  $X^j_{js}$, $X_{j\iota s}$ give rise to families $\X^{j'}_{js}$, $j\in\overline{\J}_j$,  $\X^j_{js}$, $\X_{j\iota s}$ of probability distributions  on $\Omega$; these families (and densities from the families) inherit colors from their indices $j$.  We claim that the resulting entities  form a legitimate setup for a sequential test. All we need in order to justify this claim is to verify that the $S$-component of our setup is as required, that is, that (a) for every $j\in \{1,...,J\}$, $\iota_{jS}=1$, whence $X_{j1S}=X_j$ and $L_S=J$, and that (b) $q,q'\in  \{1,...,L_S\}$ (recall that $L_S=J$) are $S$-close  if and only if $q$ and $q'$ are of the same color.
\begin{quote}{\small
To verify (a), note that $\ell_{jj'S}(\cdot)\equiv -1$ whenever $(j,j')\in\cR$, implying that $X^{j'}_{jS}=\emptyset$ when $j'\in\overline{\J}_j$ and $X^j_{jS}=X_j$, as claimed in (a).
To verify (b), note that as it was already mentioned, for $j$, $j'$ of different colors, the risk of the detector yielded by Theorem \ref{the1} as applied to the sets $X=X_j$, $Y=X_{j'}$, is $\exp\{\psi_{jj'}\}$, that is, this risk is $\leq \exp\{-d\}$. Invoking the already verified (a), we conclude that $\epsilon_{qq',S}\leq\exp\{-d\}$ whenever $1\leq q,q'\leq L_S=J$ and $q,q'$ are of different colors. As we have seen, $r(S)<d$, whence $\delta_S=\exp\{-r(S)\}>\exp\{-d\}$. The bottom line is that whenever $1\leq q,q'\leq J$ and $q,q'$ are of different colors, we have $\epsilon_{qq',S}<\delta_S$; this observation combines with the definition of $S$-closeness to imply that $q,q'$ are $S$-close if and only if $q,q'$ are of different colors, as claimed in (b).
}
\end{quote}
The legitimate setup we have presented induces a sequential test, let it be denoted by $\T$. We are about to analyse the properties of this test.
\subsection{Analysis}\label{sec:analysis}
Our  first observation is that for the sequential setup we have presented one has $k(s)\leq \bar{k}(s)$, $1\leq s\leq S$. To verify this claim, we need to check that when $k=\bar{k}(s)$, we have \[\Opt(k,s)<\epsilon_s={{\epsilon\over S}}.
 \]
As we have already mentioned  (see the comment after the definition (\ref{seq9}) of $\Opt(\cdot,\cdot)$), $\Opt(k,s)$ is the spectral norm of the entrywise nonnegative symmetric matrix $D^{ks}$ of the size $L_s\times L_s$ with entries not exceeding $\delta_s^k$. Since, by construction, $L_s\leq J^2$ {\sl and the diagonal of $D$ is zero}, the spectral norm of $D^{ks}$  does not exceed ${(J^2-1)}\delta_s^k$. The latter quantity indeed is $< \epsilon_s$ when $k=\bar{k}(s)$, see (\ref{seq11}).
\paragraph{Worst-case performance.} In the analysis to follow, we assume that $\epsilon\in(0,\four)$.
\par
By Proposition \ref{prop_sequential}, the sequential test $\T$ always accepts at most one of the hypotheses $H_1,...,H_I$, and the probability not to accept the true hypothesis is at most $\epsilon$; moreover, the number of observations used by $\T$ never exceeds $k(S)\leq \overline{K}:=\bar{k}(S)$.
 On the other hand, from the definition \rf{defd} of $d$ and  Corollary \ref{cor1} it follows that in order for a whatever test to decide on the hypotheses $H_1,...,H_I$ with risk $\epsilon$ via stationary repeated observations, the size of the observation sample should be {\sl at least}\footnote{Indeed, this is exactly the smallest number of observations which is necessary, according to Corollary \ref{cor1}, to separate with the risk $\leq \epsilon$ the pair of hypotheses corresponding to $(j,j')\in \cR$ for which $\psi_{jj'}=-d$.}
\begin{equation}\label{Kstarequals}
K^+=\left[{\half\ln(1/\epsilon)-\ln(2)\over\ln(1/\epsilon)}\right]{\ln(1/\epsilon)\over d}\geq {\ln(1/\epsilon)\over 4d}.
\end{equation}
As a result, {\sl unless $d$ is ``astronomically small'', $\overline{K}$ is within logarithmic factor of $K^+$}, implying quasi-optimal worst-case performance of the test $\T$. The precise statement is as follows:
\begin{proposition}\label{proploglog} Let $d>0$, $J\geq2$ and $\epsilon\in (0,\four)$ satisfy, for some $\kappa\geq1$,  the relation
\begin{equation}\label{letdsatisfy}
\ln(1/d)\leq \kappa \ln(J^2/\epsilon).
\end{equation}
Then
\begin{equation}\label{ifdsatisfies}
\overline{K}\leq \max\left[1,5{\kappa \ln (J^2/\epsilon)\over d}\right].
\end{equation}
\end{proposition}
For proof, see the appendix.

For all practical purposes we can assume that $d\geq 10^{-6}$, otherwise the {\sl lower} bound $K^+$ on the number of observations required by $(1-\epsilon)$-reliable test would be impractically large. Assuming $d\geq 10^{-6}$, (\ref{letdsatisfy}) is satisfied with $\kappa=5$ (recall that $\epsilon\leq \four$ and $J\geq2$). Thus, for all practical purposes we may treat the quantity $\kappa$ from the premise of Proposition \ref{proploglog} as a moderate absolute constant, implying that the upper bound $\overline{K}$ on the worst-case observation time of our $(1-\epsilon)$-reliable sequential test $\T$ indeed is within a logarithmic factor $O(1){\ln(J/\epsilon)\over \ln(1/\epsilon)}$ of the lower bound $K^+$ on the worst-case observation time of an ``ideal'' $(1-\epsilon)$-reliable test.
\begin{remark}\label{remloglog} It is easily seen that when $\bar{k}(s)$ grows with $s$ as rapidly as allowed by {\rm (\ref{monotonicity})}: $\bar{k}(s)=2^{s-1}$, the result completely similar to the one of Proposition \ref{proploglog} holds true in a much wider than {\rm (\ref{letdsatisfy})} range of values of $d$, specifically, in the range
$\ln(1/d)\leq CJ^2/\epsilon$, for a whatever constant $C\geq1$; and in this range, one has $\overline{K}\leq C'\max[1,{\ln(J^2/\epsilon)\over d}]$, with $C'$ depending solely on $C$.
\end{remark}
\paragraph{Actual performance.} For $\mu\in X:=\bigcup_{j=1}^J X_j$ let $s_*(\mu)$ be the smallest $s\leq S$ such that for some $j\leq J$ it holds $\mu\in X_{js}^j$ (see (\ref{partition})). Equivalently:
\begin{equation}\label{seq22}
s_*(\mu)=\min\left\{s: \exists j\leq J: \mu\in X_j \ \&\ \ell_{jj's}(\mu)\leq0\,\,\forall j'\in\overline{\J}_j.
\right\}\end{equation}
Note that $s_*(\mu)$ is well defined -- we have already seen that $X^{j'}_{jS}=\emptyset$ whenever $j'\in\overline{\J}_j$, so that $s=S$ is feasible for the right hand side problem in (\ref{seq22}).

Recall that in Section \ref{sectriskanal} we have associated with a distribution $p\in\bigcup_j \X_j$ a set of indices  $Q_s[p]$, $1\leq s\leq S$, and indices $s[p]$ and $q[p]$. Now we are in the situation where there is one-to-one correspondence $\mu\mapsto p_\mu$ between $\bigcup_j X_j$ and $\bigcup_j \X_j$; with this in mind (and allowing for slight abuse of notation), in what follows we set
$s[\mu]=s[p_\mu]$, $q[\mu]=q[p_\mu]$, $Q_s[\mu]=Q_s[p_\mu]$. In other words, by construction, for $\mu\in X:=\bigcup_jX_j$, $Q_s[\mu]$ is the set of  indices  $q$ of those of the sets $Z_{qs}$ (the latter sets are obtained by  linear ordering of the $L_s$ sets $X_{j\iota s}$)  which contain $\mu$, 
$s[\mu]$ is the smallest $s$ such that for some $q\in Q_s[\mu]$, all $q'\leq L_s$ which are $s$-close to $q$ are of the same color as the one of $q$, and the latter property is shared
by $q[\mu]\in Q_{s[\mu]}[\mu]$. 
 \begin{proposition}\label{propactual}
One has $$s[\mu]\leq s_*(\mu).$$
\end{proposition} Assume that the observations are drawn from  density $p_\mu(\cdot)$, $\mu\in\bigcup_jX_j$. By Proposition \ref{prop_sequential}, with $p_\mu$-probability $\geq1-\epsilon$ our sequential test $\T$ terminates in no more than $s[\mu]\leq s_*(\mu)$ steps and upon termination recovers correctly the color $i[\mu]$ of $\mu$ (i.e., accepts the true hypothesis $H_{i[\mu]}$, and only this hypothesis). Thus, if $\mu$ is ``deeply inside'' one of the sets $X_{j}$, meaning that $s_*(\mu)$ is much smaller than $S$, our sequential test will, with reliability $1-\epsilon$, identify correctly the true hypothesis $H_{i[\mu]}$ much faster than in $S$ stages.

\subsection{Selecting the cuts} Defining the cuts $\ell_{jj's}(\cdot)$ is one of principal ``degrees of freedom'' of the just described construction.  Informally speaking, we would like these cuts to result in as small sets $X^{j'}_{js}$, $j'\in\overline{\J}_j$ as possible,  thus increasing chances for the probability density underlying our observations to belong to $X^j_{js}$ for a small value of $s$ (such a value, as we remember, with probability $1-\epsilon$ upper-bounds the number of stages before termination). Thus, to improve over the basic option represented by the default cuts one may look, for instance, for cuts which minimize, given $s$ and $(j,j')\in \cR$, the $m_j$-dimensional volume of the set $X^{j'}_{js}$, where $m_j$ is the dimension of $X_j$. Setting $m=m_j$, $X=X_j$ and $Y=\{x\in X: \psi_{j'}(x)\geq-r(s)\}$, this goal boils down to solving the following geometric problem:
\begin{quote}
(*): {\sl Given a convex compact set $X\subset\bbr^m$ with a nonempty interior  and a nonempty convex compact subset $Y$ of $X$, find an affine function $\ell(x)$ such that the  linear inequality
$\ell(x)\geq0$ is valid on $Y$ and minimizes, under this requirement, the average linear size
$$
\Size(X_{\ell(\cdot)}):=\left[\mes_m(X_{\ell(\cdot)})\right]^{1/m}
$$
of the set $X_{\ell(\cdot)}=\{x\in X: \ell(x)\geq 0\}$.}
\end{quote}
Problem (*) seems to be heavily computationally intractable. We are about to present a crude {\sl sub}optimal solution to (*).
\par
We can assume w.l.o.g. that $Y$ intersects with the interior of $X$.\footnote{Indeed, otherwise we can specify $\ell(\cdot)$ as an affine function separating $Y$ and $X$, so that $\ell(x)\geq0$ when $x\in Y$ and $\ell(x)\leq0$ on $X$ and  $\ell(\cdot)$ is nonconstant on $X$. Clearly, $\ell(\cdot)$ is a feasible solution to (*) with $\Size(X_{\ell(\cdot)})=0$.} Let us equip $X$ with a $\vartheta$-self-concordant barrier $F(\cdot)$.\footnote{That is, $F$ is a three times continuously differentiable strictly convex function on $\inter X$ which is an interior penalty for $X$ (i.e., diverges  to $+\infty$ along every sequence of interior points of $X$ converging to a boundary point of $X$), and, in addition, satisfies specific differential inequalities; for precise definition and related facts to be used in the sequel, see \cite{NN1994}. A reader will not lose much when assuming that $X$ is a convex compact set given by a strictly feasible system $f_i(x)\leq b_i$, $1\leq i\leq m$,  of convex quadratic inequalities, $F(x)=-\sum_{i=1}^m\ln(b_i-f_i(x))$ and $\vartheta=m$.} Then the minimizer $\bar{x}$ of $F$ on $Y$ is uniquely defined, belongs to $\inter X$  and can be found efficiently by solving the (solvable) convex optimization problem $\min_{x\in Y}F(x)$. By optimality conditions for the latter problem, the affine function
$\bar{\ell}(x)=\langle \nabla F(\bar{x}),x-\bar{x}\rangle$ is nonnegative on $Y$ and thus is a feasible solution to (*). Further, from the basic facts of the theory of self-concordant barriers it follows that
\begin{itemize}
\item [A.] The {\sl Dikin ellipsoid} of $F$ at $\bar{x}$ -- the set $D=\{x\in X:\, \langle x-\bar{x},\nabla^2F(\bar{x})[x-\bar{x}]\rangle \leq 1\}$ -- is contained in $X$;
\item [B.] The set $X_{\bar{\ell}(\cdot)}=\{x\in X: \,\langle x-\bar{x},\nabla F(\bar{x})\rangle \geq0\}$ is contained in the set\footnote{What is $\vartheta+2\sqrt{\vartheta}$ below is $3\vartheta$ in \cite{NN1994}; refinement $3\vartheta\to\vartheta+2\sqrt{\vartheta}$ is due to F. Jarre, see \cite[Lemma 3.2.1]{NLN}.}
$$D^+=\{x\in\bbr^m:\,  \langle x-\bar{x},\nabla^2F(\bar{x})[x-\bar{x}]\rangle \leq \rho^2:=(\vartheta+2\sqrt{\vartheta})^2,\;\langle x-\bar{x},\nabla F(\bar{x})\rangle\geq0\}.$$
\end{itemize}
From B it follows that
\begin{equation}\label{volumeissmall}
\Size(X_{\bar{\ell}(\cdot)})\leq \Size(D^+)=\rho\Size(D'),\,\,D'=\{x\in D: \langle x-\bar{x},\nabla F(\bar{x})\rangle\geq0\}.
\end{equation}
On the other hand, if $\ell(\cdot)$ is a feasible solution to (*), then the set $X_{\ell(\cdot)}$ contains $Y$ and thus contains $\bar{x}$, implying that
$\ell(\bar{x})\geq0$. Consequently, by A we have
$$
\bar{D}:=\{x\in D:\ell(x)\geq0\}\subset \{x\in X:\ell(x)\geq0\}=X_{\ell(\cdot)},
$$
whence
$$
\Size(X_{\ell(\cdot)})\geq\Size(\bar{D}).$$
         Since $\ell(\bar{x})\geq0$ and $\ell(\cdot)$ is affine, we have $\bar{D}\supset D'':=\{x\in D:\ell(x)\geq\ell(\bar{x})\}$, and the $m$-dimensional volume of $D''$
         is at least half of the $m$-dimensional volume of $D'$ (since every one of the sets $D''$, $D'$ is either the entire ellipsoid $D$, or is the intersection of $D$ with half-space with the boundary hyperplane passing through the center of $D$).  It follows that
         \[\Size(\bar{D})\geq\Size(D'')\geq 2^{-1/m}\Size(D').
         \]
 Thus, for every feasible solution $\ell(\cdot)$ to (*) it holds
$$
\Size(X_{\ell(\cdot)})\geq \Size(\bar{D})\geq 2^{-1/m}\Size(D')\geq 2^{-1/m}\rho^{-1} \Size(X_{\bar{\ell}(\cdot)}),
$$
where the concluding $\geq$ is due to (\ref{volumeissmall}). We see that the feasible solution $\bar{\ell}(\cdot)$ to (*) (which can be found efficiently) is optimal within the factor $2^{1/m}\rho=2^{1/m}[\vartheta+2\sqrt{\vartheta}]$. This factor is moderate when $X$ is an ellipsoid (or the intersection of $O(1)$ ellipsoids), and can be unpleasantly large when $\vartheta$ is large. However, this is the best known to us computationally tractable approximation to the optimal solution to (*).
\subsection{Application in the Gaussian case}
\begin{figure}
$$
\epsfxsize=180pt\epsfysize=170pt\epsffile{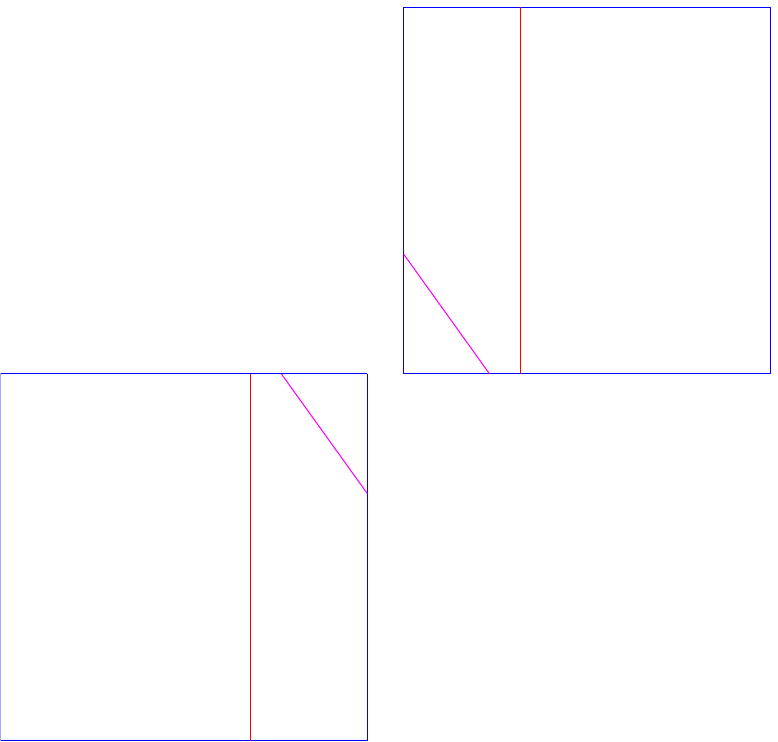}
$$
\caption{\label{newfig} \small $X_1$, $X_2$ (blue squares) and cuts (red -- default, magenta -- smart) for $s=11$ \\
{[$n=2$, $\epsilon=0.01$, $\delta=0.1$, $S=14$, $r(11)=0.0092$]}.}
\end{figure}

{
\paragraph{Comparing cut policies.} To get impression on the effect of different cut policies (default cuts vs. cuts based on self-concordant barriers), consider the Gaussian case with $J=2$, $X_1=\{x\in\bbr^n:\delta\leq x_1\leq 1+\delta,0\leq x_i\leq1,\,2\leq i\leq n\}$, $X_2=\{x\in\bbr^n:-1\leq x_i\leq 0,\,1\leq i\leq n\}$, $I_1=\{1\}$, $J_2=\{2\}$ (i.e., color of $X_1$ is 1, and color of $X_2$ is 2), and let $\bar{k}(s)=2^{s-1}$, $s=1,2,...$.  We have
$$
\begin{array}{l}
\psi(\mu,\nu)=-{1\over 8}\|\mu-\nu\|_2^2,\,\,
\psi_1(\nu)=-{1\over 8}\left[(\nu_1-\delta)^2+\sum_{i=2}^n\nu_i^2\right],\,\,\psi_2(\mu)=-{1\over 8}\|\mu\|_2^2;\\
\psi_{1,2}=\psi_{2,1}=-{1\over 8}\delta^2,\,\,d:=\min[-\psi_{1,2},-\psi_{2,1}]={\delta^2\over 8};\\
a_{1,2}=b_{2,1}=e:=[1;0;...;0],\,\,b_{1,2}=a_{2,1}=0\in\bbr^n,\,\,
e_{1,2}=f_{2,1}=e,\,\,f_{1,2}=e_{2,1}=0.\\
\end{array}
$$
Further,
$$
r(s)={2\ln(4S/\epsilon)\over 2^s},\,1\leq s\leq S,
$$ where $\epsilon$ is the target risk, and $S$ is the smallest positive integer such that
$r(S)$ as given by the above formula is $<\chi={\delta^2\over 8}$
\par
Now, the default cuts are (see (\ref{implyingcut}))
$$
\begin{array}{ll}
(1,2,r):&\ell_{1,2,r}(\mu):=-\mu_1+2\sqrt{2r}\leq0\\
(2,1,r):&\ell_{2,1,r}(\nu):=-\delta+\nu_1+2\sqrt{2r}\leq0\\
\end{array}
$$
and the sets $Z_{qs}$ associated with the default cuts are:
$$
\begin{array}{ll}
\hbox{color \# 1:}&\left\{\begin{array}{rcl}Z_{1s}&=&X_{1,s}^1=\{\mu:\max[\delta,2\sqrt{2r(s)}]\leq\mu_1\leq 1+\delta\ \&\ 0\leq \mu_i\leq1,\,2\leq i\leq n\},\\
Z_{2s}&=&X_{1,s}^2=\{\mu:\delta\leq \mu_1\leq 2\sqrt{2r(s)}\ \& \ 0\leq mu_i\leq 1,\,2\leq i\leq n\};\\
\end{array}\right.\\
\hbox{color \# 2:}&\left\{\begin{array}{rcl}Z_{3s}&=&X_{2,s}^2=\{\nu:-1\leq \nu_1\leq \min[0,\delta-2\sqrt{2r(s)}]\ \&\ -1\leq \nu_i\leq 0,\,2\leq i\leq n\},\\
Z_{4s}&=&X_{2,s}^1=\{\nu:\delta-2\sqrt{2r(s)}\leq \nu_1\leq0\ \&\ -1\leq \nu_i\leq 0,\,2\leq i\leq n\}.\\
\end{array}\right..\\
\end{array}
$$
The ``good'' sets here are $Z_{1s}$ and $Z_{3s}$, meaning that if a single observation is distributed  according to $\cN(\mu,I_d)$ with $\mu$ belonging to $Z_{1s}$ or $Z_{3s}$, our sequential test $\T$ with probability at least $1-\epsilon$ terminates  in course of the first $s$ stages with correct conclusion on the color of $\mu$. We see also that when $\delta\ll 1$ and $s$ is ``moderate,'' meaning that $r(s)\gg\delta^2$, the bad sets $Z_{2s}$, $Z_{4s}$, in terms of their $n$-dimensional volume, form $O(1)\sqrt{r(s)}$-fractions of the respective boxes $X_1$, $X_2$, provided default cuts are used. When using ``smart'' cuts -- those induced by self-concordant barriers for $X_1$, $X_2$ -- it is immediately seen that in the range $1\gg\delta\gg r(s)\gg\delta^2$ (what exactly $\gg$ means, depends on $n$), the bad sets $Z_{2s}$ and $Z_{4s}$ become simplexes:
$$
\begin{array}{rcl}
Z_{2s}&=&\{\mu=[\delta;0;...;0]+[\lambda_1;...;\lambda_n]:\lambda\geq0,\sum_{i=1}^nc_i\lambda_i\leq \sqrt{r(s)n}\},\\
Z_{4s}&=&\{\nu=-[\lambda_1;...;\lambda_n]:\lambda\geq0,\sum_{i=1}^nc_i\lambda_i\leq \sqrt{r(s)n}\},\\
\end{array}
$$
where $c_i=c_i(n,s)$ are of order of 1; the good sets $Z_{1s}$, $Z_{3s}$ are the closures of the complements of these simplexes to the respective boxes $X_1$, $X_2$. The new bad sets are much smaller than the old ones; e.g., their volumes are of order of $r^{n/2}(s)$ -- much smaller than the volumes $O(\sqrt{r(s)})$ of the old bad sets, see Figure \ref{newfig} and Table \ref{newtable}.
}

\begin{table}
\begin{center}
\begin{tabular}{|c|c|c|c|c|c|c|}
\hline
cuts&$n=2$&$n=3$&$n=4$&$n=5$&$n=6$&$n=20$\\
\hline
\hline
default&0.32&0.32&0.32&0.32&0.32&0.32\\
\hline
smart&3.8e-2&7.6e-3&1.4e-3&2.2e-4&3.4e-5&$<$1.3e-18\\
\hline
\end{tabular}
\end{center}
\caption{\label{newtable}\small $n$-dimensional volumes of ``bad'' sets $Z_{qs}$ for $s=11$
{[$\epsilon=0.01$, $\delta=0.1$, $S=14$, $r(11)=0.0092$]}.}
\end{table}

\paragraph{Upper bounding $s[\mu]$.}
{
In Gaussian case {\sl with default cuts}, the quantity $s_*(\mu)$ which, as we have just seen, is a $(1-\epsilon)$-reliable upper bound on the number of stages in which $\T$ recognizes $(1-\epsilon)$-reliably  the true hypothesis $\H_{i[\mu]}$ provided the observations are drawn from  $p_\mu(\cdot)$, admits a transparent geometric upper bound. Observe, first, that in Gaussian case we have
\begin{equation}\label{gaussianchi}
\chi={1\over 8}\min_{(j,j')\in \cR}\min_{a\in X_j,b\in X_{j'}}\|a-b\|_2^2.
\end{equation}
Now let $\rho(\mu)$ be the largest $\rho$ such the $\|\cdot\|_2$-ball of radius $\rho$ centered at $\mu$ is contained in certain $X_j$. We claim that
\begin{equation}\label{gauss222}
s[\mu]\leq s_*(\mu)\leq \bar{s}(\mu):=\min\{s\leq S: 2\sqrt{2r(s)}\leq 2\sqrt{2\chi}+\rho(\mu)\},
\end{equation}
meaning that the deeper $\mu$ is ``inside'' one of $X_j$ (the larger is $\rho(\mu)$), the smaller is the number of observations needed for $\T$ to identify correctly the color of $\mu$.
\begin{quote}
{\small
Justification of (\ref{gauss222}) is as follows. The first inequality in (\ref{gauss222}) has already been proved. To prove the second inequality, observe, first, that $\bar{s}:=\bar{s}(\mu)$ is well defined (indeed, as we have seen, $r(S)<\chi)$. Let $j$ be such that the $\|\cdot\|_2$-ball $B$ of radius $\rho(\mu)$ centered at $\mu$ is contained in $X_{j}$, and let $j'\in\J^o_{j}$. In the notation from the description of the Gaussian case default cuts we have, by (\ref{whatisetacut}),  $\eta_{jj'}^T[\mu'-a_{j_*j}]\geq0$
for all $\mu'\in X_{j}$ and thus for all $\mu'\in B$, and therefore $\eta_{jj'}^T[\mu-a_{jj'}]\geq\rho(\mu)$, since by construction $\|\eta_{jj'}\|_2=1$. Further, $\delta_{jj'}$ is the $\|\cdot\|_2$-distance between $X_{j}$ and $X_{j'}$, whence $\delta_{jj'}\geq2\sqrt{2\chi}$ due to (\ref{gaussianchi}) combined with $j'\in\C^o_{j}$. Thus,
\begin{equation}\label{seq401}
\forall (j'\in \J^o_{j}): \delta_{jj'}+\eta_{jj'}^T(\mu-a_{jj'})\geq 2\sqrt{2\chi}+\rho(\mu).
\end{equation}
Invoking the definition of $\bar{s}:=\bar{s}(\mu)$ and the description (\ref{implyingcut}) of $(jj'r(\bar{s}))$-cuts, we conclude from (\ref{seq401}) that $\mu$ indeed satisfies all these cuts and therefore $\mu\in X_{j\bar{s}}^j$, that is, $s_*(\mu)\leq\bar{s}$ (recall the definition of $s_*(\mu)$), as claimed.}
\end{quote}
}

\begin{figure}
$$
\begin{array}{cc}
\begin{array}{c}
\\[-30pt]
\epsfxsize=140pt\epsfysize=100pt\epsffile{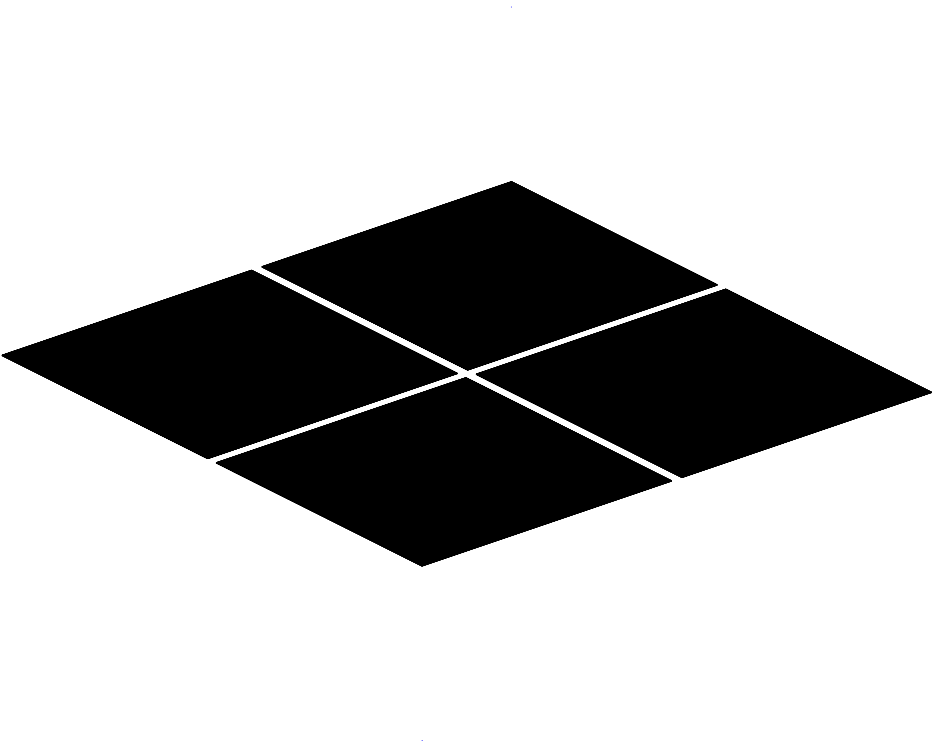}\\
\end{array}&\epsfxsize=230pt\epsfysize=140pt\epsffile{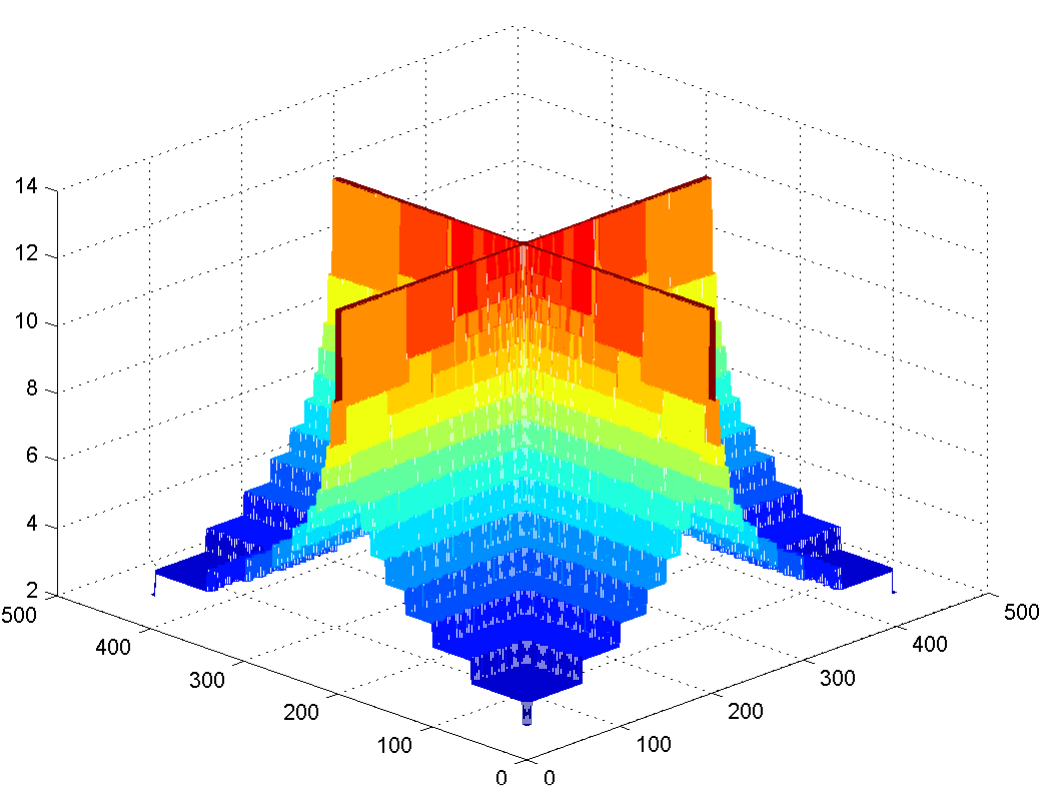}\\
\end{array}
$$
\caption{\label{figseqprof} ``Power'' of sequential test, Gaussian case. Left display: sets $X_1,X_2,X_3,X_4$;  right display:  $\ln (k(s_*(\mu)))$  as a function of $\mu\in\bigcup\limits_{i=1}^4X_i$, $\epsilon=0.01$.}
\end{figure}

\paragraph{Numerical illustration.} The following numerical experiment highlights the power of sequential testing in the Gaussian case with default cuts. In this experiment, we are given $J=4$ sets $X_j\subset\bbr^2$; $X_1$ is the square $\{0.01\leq x_1,x_2\leq 1\}$, $X_2,X_3,X_4$ are obtained from $X_1$ by reflections w.r.t. the coordinate axes and the origin. The partition of the index set into groups $\J_i$ is trivial -- these groups are the elements of $\J=\{1,2,3,4\}$, so that our goal is to recognize which of the sets $X_j$ contains the mean $\mu$ of the observation. Figure \ref{figseqprof} presents the graph of the {\sl logarithm} of the $0.99$-reliable upper bound on the number of observations used by the sequential test with $S=20$, $\bar{k}(s)=2^{s-1}$ and  $d=5.0$e$-5$ as a function of the mean $\mu$ of the observation. We see that the savings from sequential testing are quite significant. The related numbers are as follows: when selecting $\mu$ in $\bigcup\limits_{j=1}^4 X_j$ at random according to the uniform distribution, the empirical average of the number of observations before termination is as large as $1.6\cdot10^5$, reflecting the fact that $X_j$ are pretty close to each other. At the same time, the median number of observations before termination is just 154, reflecting the fact that in our experiment $\mu$, with reasonably high probability, indeed is ``deeply inside'' the set $X_j$ containing $\mu$.

\appendix
\section{Proof of Proposition \ref{prop_sequential}}
Let $p_*$, $j_*$, $i_*$ be the entities from the proposition. For $1\leq s\leq S$, there exist $q_s\in\{1,...,L_s\}$ such that $p_*\in \Z_{q_ss}$; without loss of generality, we can assume that $q_{s[p_*]}=q[p_*]$. For $1\leq s\leq S$ let
$$
\begin{array}{rcl}
\E_1^s&=&\{\omega^K:\, T_{s} \hbox{\ as applied to $\omega^{k(s)}$ does not accept the hypothesis $H_{q_ss}$}\},\\
\E_2^s&=&\{\omega^K: \,T_s \hbox{\ as applied to $\omega^{k(s)}$ accepts a hypothesis $H_{qs}$ with $q$ not $\C^s$-close to $q_s$}\},\\
{ \E^s}&=&\E_1^s\cup E_2^s.\\
\end{array}
$$
From the just outlined properties of $T_s$ it follows that $p_*$-probability { of $\E^s$ does} not exceed $\epsilon_s$. Now let
$$
\E^*=\{\omega^K: \hbox{\ no one of the events $\E_1^s$, $\E_2^s$, $1\leq s\leq S$, takes place.}\},
$$
so that $p_*$-probability of $\E^*$ is at least $1-2\sum_{s=1}^S\epsilon_s = 1-\epsilon$. All we need is to verify that $\E\supset \E^*$, which is immediate. Indeed, let $s_*=s[p_*]$, $q_*=q[p_*]$ (so that $q_*=q_{s_*}$ due to $q_{s[p_*]}=q[p_*]$), and let $\omega^K\in \E^*$. By the latter inclusion, $\E_1^{s_*}$ does not take place, what implies that \\
\indent (a) {\em $T_{s_*}$ as applied to $\omega^{k(s_*)}$ accepts
$H_{q_*s_*}$.}\\
Besides this, since $\omega^K\in \E^*$, we have $\omega^{k(s_*)}\not\in \E_2^{s_*}$, implying, by definition of $\E_2^{s_*}$, that all hypotheses $H_{qs_*}$ accepted by $T_{s_*}$ as applied to $\omega^{k(s_*)}$ are such that $q$ is $\C^{s_*}$-close to $q_*=q[p_*]$. By definition of $s_*=s[p_*]$ and $q_*=q[p_*]$, for all $q$ which are $\C^{s_*}$-close to $q_*$, the hypotheses $H_{qs_*}$ are of the same color as the accepted  by $T_{s_*}$, by (a), hypothesis $H_{q_*s_*}$. Thus, when $\omega^K\in\E^*$, we have
\\
\indent (b) {\em the set of hypotheses $H_{qs_*}$ accepted by $T_{s_*}$ is nonempty, and all these hypotheses are of the same color, equal to the color of $\mu$.}
\\
Invoking the description of $\T$, we conclude that when $\omega^K\in \E^*$, the test $\T$ terminates not later than at step $s_*=s[p_*]$ and the termination is productive -- some of the  hypotheses $\H_i$ indeed are accepted (in fact, in this case exactly one of the hypotheses $\H_i$ is accepted, since, as it was already mentioned, $\T$ never accepts more than one hypothesis).
\par
Now let $\omega^K\in \E^*$, and let $\bar{s}$ be the step at which $\T$ terminates. As we have already seen, $\bar{s}\leq s_*=s[p_*]$ and $\T$ terminates accepting exactly one hypothesis $H_{\bar{i}}$, where $\bar{i}$ is a deterministic function of $\omega^{k(\bar{s})}$. Since $p_*$ obeys hypothesis $H_{q_{\bar{s}}\bar{s}}$ (by definition of $q_{\bar{s}}$) and $\omega^K\not\in \E_1^{\bar{s}}$, the test $T_{\bar{s}}$ as applied to $\omega^{k(\bar{s})}$ accepts the hypothesis $H_{q_{\bar{s}}\bar{s}}$ (which, by construction, has the same color $i_*$ as $p_*$). The latter observation combines with the termination rule for $\T$ to imply that the outcome of $\T$ in the case of $\omega^K\in\E^*$ is the hypothesis $H_i$ of the same color $i$ as $H_{q_{\bar{s}}\bar{s}}$, that is, the hypothesis $\_i$ with the same color as the color $i_*$ of $p_*$, that is, the true hypothesis $H_{i_*}$, and that the outcome is obtained not later than at the stage $s[p_*]$. \qed
\section{Proof of Proposition \ref{proploglog}}
Let $d>0$, $J\geq2$, $\epsilon\in(0,\four)$, $\varkappa>2$, and $\kappa\geq1$ satisfy (\ref{letdsatisfy}). Relation (\ref{ifdsatisfies}) is trivially true when $S=1$ and thus $\overline{K}=\bar{k}(1)=1$. We denote $\vartheta=J^2/\epsilon$, and
assume from now on that $S>1$, so that (\ref{seq10}) is not satisfied when $S=\bar{k}(S)=1$ and therefore $d\leq \ln(\vartheta)$.
\paragraph{1$^0$.} We claim that
\begin{equation}\label{eq12345}
S\leq \bar{S}:=\left\rfloor \sigma={4\kappa\ln(\vartheta)\over 3d}\right\lfloor.
\end{equation}
To justify this claim, it suffices to verify that
\begin{equation}\label{eq12346}
{\ln(\vartheta \bar{S})/\bar{k}(\bar{S})}<d.
\end{equation}
Indeed, we have $\sigma\geq 4/3$, whence $\sigma\leq \bar{S}\leq {3\over 2}\sigma$. Using the bound $\bar{k}(S)\geq 2\bar{S}\geq 2\sigma$, we conclude that the left hand side in (\ref{eq12346}) does not exceed $\C d$, with $\C$ given by
\[
\C={\ln\left(2\kappa\vartheta\ln(\vartheta)/d\right)\over {8\over 3}\kappa \ln(\vartheta)}\leq
{\ln\left(2\kappa\vartheta\ln(\vartheta)\right)+\kappa \ln(\vartheta)\over {8\over 3}\kappa\ln(\vartheta)},
\]
where the concluding $\leq$ is due to $\ln (1/d)\leq\kappa \ln(\vartheta)$ (the latter is nothing but \rf{letdsatisfy}). From the relations $\vartheta=J^2/\epsilon>16$ and $\kappa\geq1$ it immediately follows that
$\C<1$, implying the validity of  (\ref{eq12346}) and thus -- the validity of (\ref{eq12345}).
\paragraph{2$^0$.} Since $1<S\leq \bar{S}$, we have
\[
d\leq {\ln([S-1]\vartheta)\over \bar{k}(S-1)}\leq {\ln([\bar{S}-1]\vartheta)\over \bar{k}(S-1)}\leq {\ln(\sigma\vartheta)\over \bar{k}(S-1)}.
\] Thus
\bse
\overline{K}&=&\bar{k}(S)\leq 2\bar{k}(S-1)\leq {2\ln(\sigma\vartheta)\over d}= {2\ln({4\over 3}\kappa\vartheta\ln(\vartheta)/d)\over d}\\
\mbox{[by \rf{letdsatisfy}]}&\leq&
{2\left[\ln({4\over 3}\kappa\vartheta\ln(\vartheta))+\kappa\ln(\vartheta)\right]\over d}=\C'{\kappa\ln(\vartheta)\over d},\\
\C'&=&{2[\ln({4\over 3}\kappa\ln(\vartheta))+(\kappa+1)\ln(\vartheta)]\over\kappa\ln(\vartheta)}.
\ese
Since $\vartheta>16$ and $\kappa\geq 1$ we have
$\C'\leq 5$, and we arrive at (\ref{ifdsatisfies}).
\qed
\section{Proof of Proposition \ref{propactual}}
Let $s_*=s_*(\mu)$, and let $j_*$ be such that $\mu\in X_{j_*s_*}^{j_*}$, see the definition of $s_*(\mu)$.
The set $X_{j_*s_*}^{j_*}\ni\mu$ is some $Z_{q_*s_*}$ with $q_*\in Q_{s_*}[\mu]$. We claim that\\[3pt]
 $\quad$ (!) {\sl all $q'$ which are $s_*$-close to $q_*$ are of the same color as $q_*$};\\[3pt]
note that the validity of (!) means that setting $s=s_*$, $q=q_*$, we meet the requirements in (\ref{mayhappen}), implying, by definition of $s[\mu]$, our target relation $s[\mu]\leq s_*$.\par
To verify (!), let $q'$ be $s_*$-close to $q_*$, so that by definition of $s$-closeness either\\
 $\quad$(a) $q'$ and $q_*$ are of the same color, or\\
  $\quad$(b) $q'$ and $q_*$ are of different colors and $\epsilon_{q'q_*,s_*}>\delta_{s_*}=\exp\{-r(s_*)\}$;\\
  all we need to prove is that (b) in fact is impossible. Assume, on the contrary, that (b) takes place, and let us lead this assumption to a contradiction. As it was already mentioned, \[\epsilon_{q'q_*,s_*}=\exp\{\lambda\} \;\mbox{where}\;
  \lambda=\max_{\mu\in Z_{q_*s_*},\nu\in Z_{q's_*}}\psi(\mu,\nu).
  \]
  Thus, we are in the case when the color $i$ of $q'$
 differs from the color $i_*$ of $q_*$ and
 \begin{equation}\label{contradict}
\lambda:=\max_{\mu\in Z_{q_*s_*},\nu\in Z_{q's_*}}\psi(\mu,\nu)>-r(s_*).
 \end{equation}
 Since $q'$ is of color $i$, we have $Z_{q's_*}\subset X_{j'}$ for some $j'\in\J_i$, and since $i\neq i_*$, $j'$ and $j_*$ are of different colors, or, equivalently, $j'\in\overline{\J}_{j_*}$. Now, by the definition of $q_*$ we have
 \[
Z_{q_*s_*}=X^{j_*}_{j_*s_*}=\{\mu\in X_{j_*}: \ell_{j_*js_*}(\mu)\leq 0\,\,\forall j\in\overline{\J}_{j_*}\}.
\]
Hence, taking into account that $j'\in \overline{\J}_{j_*}$ and that $\ell_{j_*j's_*}(\cdot)\leq0$ is $(j_*j's_*)$-cut, $\mu\in Z_{q_*s_*}$ implies that
\[
\psi_{j'}(\mu)\leq -r(s_*).
\]
In its turn, the latter relation, due to $Z_{q's_*}\subset X_{j'}$ and that $\psi_{j'}(\cdot)=\max_{\nu\in X_{j'}}\psi(\mu,\nu)$, implies that
$$
\max_{\mu\in Z_{q_*s_*},\nu\in Z_{q's_*}}\psi(\mu,\nu)\leq-r(s_*),
$$
which contradicts (\ref{contradict}). 
\qed
\end{document}